\documentclass[11pt,a4paper]{amsart}

\usepackage{vmargin}
\setmargrb{2cm}{2cm}{2cm}{2cm}

\usepackage{amsmath, amssymb, amsfonts, amsthm}
\usepackage{txfonts}
\usepackage{url}
\usepackage{graphicx}
\usepackage[abs]{overpic}
\usepackage{subfigure}

\theoremstyle{plain}
\newtheorem{thm}{Theorem}[section]

\newtheorem{prop}[thm]{Proposition}
\newtheorem{cor}[thm]{Corollary}

\theoremstyle{definition}
\newtheorem{defn}[thm]{Definition}
\newtheorem{exmp}[thm]{Example}

\theoremstyle{remark}
\newtheorem{rem}[thm]{Remark}

\newcommand{\FF}{\mathbb{F}}
\newcommand{\RR}{\mathbb{R}}
\newcommand{\ZZ}{\mathbb{Z}}
\newcommand{\conv}{\operatorname{conv}}
\newcommand{\Aut}{\operatorname{Aut}}
\newcommand{\link}{\operatorname{link}}
\renewcommand{\star}{\operatorname{star}}
\newcommand{\del}{\operatorname{del}}
\newcommand{\dw}{\operatorname{DW}}
\newcommand{\susp}{\operatorname{Susp}}
\newcommand\SetOf[2]{\left\{#1\,\vphantom{#2}\right|\left.\vphantom{#1}\,#2\right\}}

\begin{document}

\title{One-Point Suspensions and Wreath Products of Polytopes and Spheres}
\author[Joswig and Lutz]{Michael Joswig \and Frank H. Lutz}
\address{Michael Joswig, Institut f\"ur Mathematik, MA 6-2, TU Berlin,
  10623 Berlin, Germany}
\address{Frank H. Lutz, Institut f\"ur Mathematik, MA 6-2, TU Berlin,
  10623 Berlin, Germany}

\begin{abstract}
  It is known that the suspension of a simplicial complex can be realized with only one additional point.  Suitable
  iterations of this construction generate highly symmetric simplicial complexes with various interesting combinatorial
  and topological properties.  In particular, infinitely many non-PL spheres as well as
  contractible simplicial complexes with a vertex-transitive group of automorphisms 
  can be obtained in this way.
\end{abstract}

\maketitle

\section{Introduction}
McMullen \cite{McMullen1976} constructed projectively unique convex polytopes as the joint convex hulls of polytopes in
mutually skew affine subspaces which are attached to the vertices of yet another polytope.  It is immediate that if the
polytopes attached are pairwise isomorphic one can obtain polytopes with a large group of automorphisms.  In fact, if
the polytopes attached are simplices, then the resulting polytope can be obtained by successive wedging (or rather its
dual operation).  This dual wedge, first introduced and exploited by Adler and Dantzig \cite{AdlerDantzig1974} in 1974
for the study of the Hirsch conjecture of linear programming, is essentially the same as the one-point suspension in
combinatorial topology.  It is striking that this simple construction makes several appearances in the literature, while
it seems that never before it had been the focus of research for its own sake.  The purpose of this paper is to collect
what is known (for the polytopal as well as the combinatorial constructions) and to fill in several gaps, most notably
by introducing \emph{wreath products} of simplicial complexes.

In particular, we give a detailed analysis of wreath products in order to provide explicit descriptions of highly
symmetric polytopes which previously had been implicit in McMullen's construction. This is instrumental in proving that
certain simplicial spheres that occurred in the process of enumerating the types of combinatorial manifolds with few
vertices are, in fact, polytopal.

Non-PL spheres have been constructed by Edwards~\cite{Edwards1975} and Cannon~\cite{Cannon1979} by suspending (at least
twice) any arbitrary homology sphere.  By enumeration, Lutz \cite{Lutz2003pre} obtained three $17$-vertex triangulations
of the Poincar\'e homology $3$-sphere with a vertex-transitive group action.  The wreath products of the boundary of a
simplex with these triangulations form a new class of non-PL-spheres with a vertex-transitive automorphism group.

It is a --- presumably difficult --- open problem to decide whether or not there exist vertex-transitive non-evasive
simplicial complexes.  It is even unclear if vertex-transitive collapsible complexes exist. If not, then this would
settle the long-standing evasiveness conjecture for graph properties of complexity theory; see Kahn, Saks, and
Sturtevant \cite{KahnSaksSturtevant1984}.  Few vertex-transitive contractible and $\ZZ$-acyclic complexes are known.  A
new family of vertex-transitive contractible simplicial complexes arises via the wreath product construction.  However,
we can show that a non-evasive wreath product necessarily has a non-evasive factor.  Thus wreath products do not lead to
a solution of the evasiveness conjecture.

One-point suspensions have recently been employed successfully to construct non-constructible, non-shell\-able, not
vertex-decomposable, as well as non-PL spheres with few vertices; see \cite{BjoernerLutz2000}, \cite{BjoernerLutz2003},
\cite{Lutz2003gpre}, \cite{Lutz2004a}, and \cite{Lutz2004b}.  Here, we will investigate, how these combinatorial
properties are respected by one-point suspensions and wreath products.

\section{The Polytopal Constructions}\label{sec:polytopal}
A \emph{convex polytope} is the convex hull of finitely many points in $\RR^d$ (interior description) or, equivalently,
the bounded intersection of finitely many affine halfspaces (exterior description).  The two descriptions are dual to
each other by means of cone polarity.  The \emph{dimension} of a polytope is the dimension of its affine span.  A
\emph{vertex} is a point of a polytope which is not redundant in its interior description.  Dually, for a
full-dimensional polytope~$P$, a \emph{facet} is the intersection of $P$ with the boundary hyperplane of an affine
halfspace which is not redundant in the exterior description of~$P$.  For an introduction to polytope theory the reader
is referred to Ziegler~\cite{Ziegler1995}.

\subsection{The dual wedge of a polytope}\label{sec:dual-wedge}
Let $P\subset\RR^d$ be a $d$-dimensional polytope (or $d$-polytope for short), and let $v$ be a vertex.  The
$(d+1)$-polytope \[\dw(v,P)=\conv(P\oplus 0\cup \{v\oplus 1,v\oplus(-1)\})\] is called the \emph{dual wedge} of $P$ with
respect to~$v$.  It has the same vertices as $P$ (embedded into $\RR^{d+1}$: the notation ``$\oplus$'' is used to
indicate which additional coordinate to append), except for $v$ which splits into an ``upper'' copy $v\oplus 1$ and a
``lower'' copy $v\oplus(-1)$.  The facets of $\dw(v,P)$ are the following: For each facet $F$ of~$P$ which does not
contain~$v$ we obtain an upper cone $\conv(F\oplus 0\cup v\oplus 1)$ and a lower cone $\conv(F\oplus 0\cup
v\oplus(-1))$.  And each facet $G$ which contains~$v$ re-appears as its dual wedge $\dw(v,G)$.  Since the dual wedge of
a point clearly is a line segment, we recursively obtain a complete combinatorial description.  In particular, the dual
wedge \emph{is} a combinatorial construction: Given two polytopes $P$, $P'$ and a combinatorial isomorphism $\phi:P\to
P'$ the dual wedges $\dw(v,P)$ and $\dw(\phi(v),P')$ are combinatorially isomorphic for any vertex $v$ of~$P$.

The dual wedge of a line segment, with respect to any one of its two vertices, is a triangle.  Therefore, the recursive
description immediately implies that $\dw(v,P)$ contains a triangular $2$-face if $d\ge2$.  Moreover, $\dw(v,P)$ is a
$(d+1)$-simplex if and only if $P$ is a $d$-simplex.  This further implies that $\dw(v,P)$ is simplicial if and only if $P$
is.

The reflection at the hyperplane $x_{d+1}=0$ in~$\RR^{d+1}$ interchanges $v\oplus 1$ with $v\oplus(-1)$ and fixes all
other vertices of~$\dw(v,P)$.

Below we especially focus on iterated dual wedge constructions.

\begin{prop}
  Let $P$ be a $d$-polytope with a vertex~$v$.  Then the $(d+2)$-polytopes $\dw(v\oplus 1,\dw(v,P))$ and
  $\dw(v\oplus(-1),\dw(v,P))$ are isometric.
\end{prop}

Typically we are only interested in the combinatorial type of a dual wedge.  Hence we abbreviate $\dw^2(v,P)$ for either
$\dw(v\oplus 1,\dw(v,P))$ or $\dw(v\oplus(-1),\dw(v,P))$.  Likewise we write $\dw^k(v,P)$ for further iterations.

\subsection{The wreath product of polytopes}
Let $P\subset\RR^d$ be a $d$-polytope, and let $Q\subset\RR^e$ be an $e$-polytope.  Just in order to simplify the
description we assume that the vertex barycenters of both, $P$ and~$Q$, are zero.  Let $v_1,\dots,v_m$ be the vertices
of~$P$, and let $w_1,\dots,w_n$ be the vertices of~$Q$.  For $k\in\{1,\dots,n\}$ and $p\in\RR^d$ we define a vector
$p^k\in\RR^{nd}$ as follows: Identifying $\RR^{nd}$ with the set of matrices with $n$~rows and $d$~columns, we let $p^k$
be the $(n\times d)$-matrix with the $k$-th row equal to~$p$ and all other rows equal to zero.  Then we call the polytope
\[P\wr Q=\conv\SetOf{(v_i)^k\oplus w_k}{1\le i\le m,\ 1\le k\le n}\subset\RR^{nd+e}\]
the \emph{wreath product} of $P$ with~$Q$.  Clearly, the wreath product is full-dimensional and it has $mn$~vertices.

We use the exponent notation also for subsets of~$\RR^d$.  Moreover, we write the joint convex hull of disjoint
polytopes $R,S\subset\RR^{nd+e}$ as the \emph{join product} $R*S$.

\begin{prop}
  Take a facet $G$ of~$Q$, and assume that $w_1,\dots,w_g$ are the vertices of~$G$.  For each $k>g$, that is, for each
  vertex $w_k$ of~$Q$ which is not contained in~$G$, choose some facet $F_k$ of~$P$.  Then the iterated join
  \[F=(P^1\oplus w_1)*\dots*(P^g\oplus w_g)*((F_{g+1})^{g+1}\oplus w_{g+1})*\dots*((F_n)^n\oplus w_n)\]
  is a facet of
  $P\wr Q$, and all facets arise in this way.  We denote $F$ by $(F_{g+1},\dots,F_n;G)$.
\end{prop}

\begin{rem}
  The property that the polytopes $P$ and $Q$ both have the origin as their vertex barycenters is not strictly necessary
  in order to obtain a valid facet description as above: It suffices that the origin is an interior point.  However, the
  vertex barycenter is a fixed point of any affine transformation of a polytope, and this way, all affine
  transformations become linear.
\end{rem}

We continue with the notation of the previous proposition.  Since the vertex barycenter of~$Q$ is the origin, there is a
unique non-zero vector $\gamma\in\RR^e$ such that the linear inequality corresponding to~$G$, with indeterminate~$x$, is
$1+\langle x,\gamma\rangle\ge 0$.  Call $\gamma$ the \emph{normalized facet normal} vector of~$G$.  Similarly, let
$\phi_{g+1},\dots,\phi_n\in\RR^d$ be the normalized facet normal vectors of the facets $F_{g+1},\dots,F_n$,
respectively.  It is easy to verify that \[\left(\sum_{k=g+1}^n(1+\langle
  w_k,\gamma\rangle)(\phi_k)^k\right)\oplus\gamma\in\RR^{nd+e}\]
is the normalized facet normal vector of~$F$.

\begin{cor}
  The wreath product $P\wr Q$ is simplicial if and only if $P$ is a simplex and $Q$ is simplicial.  Moreover, $P\wr Q$
  is a simplex if and only if $P$ and $Q$ both are simplices.
\end{cor}

In general, there is no closed formula known for the $f$-vector of the wreath product.  For the important special cases
of $Q$ being either simplicial or cubical we can, however, easily count the number of facets.

\begin{cor}\label{cor:number-of-facets}
  Assume that each facet of~$Q$ has the same number of vertices, say~$c$.  Then the number of facets of $P\wr Q$ equals
  $f_{e-1}^Q(f_{d-1}^P)^{n-c}$.
\end{cor}

\begin{figure}[htbp]
  \begin{center}
    \begin{overpic}[width=6cm]{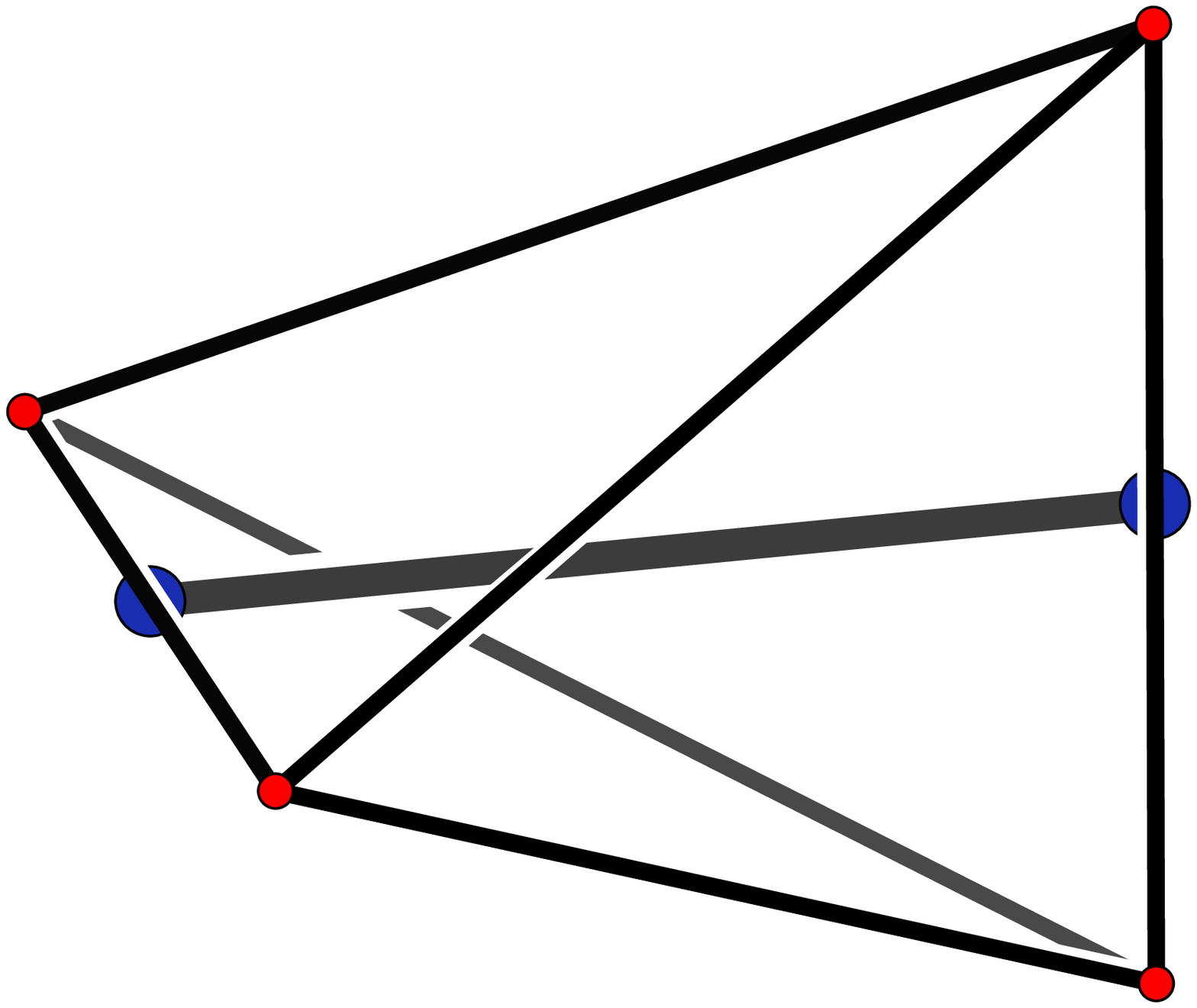}
      \put(6,28){$(-1)^1\oplus(-1)$}
      \put(-28,100){$(1)^1\oplus(-1)$}
      \put(160,10){$(-1)^2\oplus 1$}
      \put(160,147){$(1)^2\oplus 1$}
      \put(-20,60){$0\oplus(-1)$}
      \put(164,80){$0\oplus 1$}
    \end{overpic}
    \caption{Wreath product $P\wr Q$ for $P=Q=[-1,1]$.  The wreath product contains an isometric copy of
      $Q=\conv\{0\oplus(-1),0\oplus 1\}$ as shown.  For each vertex of $Q$ the boundary of the wreath product contains an
      isometric copy of~$P$: $\conv\{(-1)^1\oplus(-1),(1)^1\oplus(-1)\}$ and $\conv\{(-1)^2\oplus(-1),(1)^2\oplus 1\}$,
      respectively.}
    \label{fig:segment_wr_segment}
  \end{center}
\end{figure}

If $Q$ is a point, then $P\wr Q=P$; likewise, if $P$ is a point, then $P\wr Q=Q$.  So the first non-trivial case is
$P=Q=[-1,1]$ and
\[P\wr Q=\conv\{(-1,0,-1),(1,0,-1),(0,-1,1),(0,1,1)\}\]
is a (non-regular) tetrahedron, see Figure~\ref{fig:segment_wr_segment}.

Our terminology is justified by the following observation.

\begin{prop}\label{prop:symmetry_wreath_product_polytope}
  The wreath product of the automorphism groups $\Aut P\wr\Aut Q=(\Aut P)^n\rtimes\Aut Q$ (where the semi-direct product
  $\rtimes$ is taken with respect to the natural action of $\Aut Q$ on the $n$~vertices of~$Q$) acts as a group of
  automorphisms of $P\wr Q$.  In particular, if $\Aut P$ and $\Aut Q$ both act transitively on the set of vertices
  of~$P$ and~$Q$, respectively, then also $P\wr Q$ admits a vertex transitive group of automorphisms.
\end{prop}

The example $[-1,1]\wr[-1,1]$ above shows that the whole group of automorphisms of the wreath product can, in fact, be
larger: $(\ZZ/2)\wr(\ZZ/2)$ is the quaternion group of order eight, while the automorphism group of the $3$-simplex is
the symmetric group of degree~$4$.

One interesting map is the linear projection $\pi:P\wr Q\to P\times Q$ induced by $(v_i)^k\oplus w_k\mapsto v_i\oplus w_k$.
Additionally we define the \emph{blocking map} \[\beta:(P\wr Q)^*\to Q^*:(F_{g+1},\dots,F_n;G)^*\mapsto G^*,\] which is a
linear map between the polar polytopes.

Of special interest is the case where the first factor in the wreath product is a simplex.

\begin{prop}
  Let $\Delta_d$ be a $d$-simplex, and let $Q\subset\RR^e$ be an $e$-polytope with vertices $w_1,\dots,w_n$.  Then the
  wreath product $\Delta_d\wr Q$ is combinatorially isomorphic to the iterated dual wedge
  \[\dw^d(w_1,\dw^d(w_2,\dots \dw^d(w_n,Q)\dots)).\]
\end{prop}

\section{The Combinatorial Constructions}
Combinatorially, the dual wedge $\dw(v,P)$ of a simplicial polytope $P$ with respect to a vertex $v$ can be described as
a \emph{one-point suspension} of the boundary sphere $\partial P$ of $P$ with respect to $v$. As we will see, also the
wreath product construction $\Delta_d\wr Q$ of a $d$-dimensional simplex $\Delta_d$ with a simplicial polytope $Q$ has a
natural generalization to simplicial complexes.  For a survey on combinatorial properties of simplicial complexes see
Bj\"orner \cite{Bjoerner1995}.

\subsection{One-point Suspensions, Reduced Joins, and Wreath Products of Simplicial Complexes}
In the following, we consider finite simplicial complexes $K\neq\emptyset$.  The \emph{link}, the \emph{star}, and the
\emph{deletion} of a vertex $v$ of $K$ are the subcomplexes of $K$
\[
\begin{array}{lll}
\link_{K}(v) &:=& \SetOf{ F\in K }{ \text{$v\notin F$ and $\{v\}\cup F \in K$} },\\
\star_{K}(v) &:=& \SetOf{ F\in K }{ v\in F },\\
\del_{K}(v)  &:=& \SetOf{ F\in K }{ v\notin F },
\end{array}
\]
respectively.

\begin{defn}
  Let $K$ be a simplicial complex and let $v$ be a vertex of $K$.  The \emph{one-point suspension} $\susp_1(v,K)$ of $K$
  with respect to $v$ is the simplicial complex
  \[\susp_1(v,K):=((\partial\,\overline{v'v''}*K)\,\backslash\,(\partial\,\overline{v'v''}*\star_K(v)))\,\cup\, \overline{v'v''}*\link_K(v),\]
  where $v'$ and $v''$ are two copies of the vertex $v$ that are not contained in $K$ and which span the edge
  $\overline{v'v''}$.
\end{defn}

\begin{rem}\label{rem:susp-facets}
  The facets of $\susp_1(v,K)$ come in three kinds, depending on whether they contain $v'$, $v''$, or both: for each
  facet $F$ of~$K$ which does not contain~$v$, we obtain two coned copies $v'*F$ and $v''*F$, and for each facet $G$
  which contains~$v$, we obtain one coned copy $\overline{v'v''}*(G\setminus \{v\})$.  The canonical projection $\beta$
  which maps the facets of the one-point suspension to its base space by letting $\beta(v'*F)=\beta(v''*F)=F$ and
  $\beta(\overline{v'v''}*(G\setminus \{v\}))=G$ is \emph{not} a simplicial map; it induces a retraction of the space
  $|\susp_1(v,K)|\setminus\{v',v''\}$ to $K$, where $|\susp_1(v,K)|$ denotes a geometric realization of $\susp_1(v,K)$.
\end{rem}

Since the \emph{standard suspension} $S^0*K$ of $K$, i.e., the join product of $K$ with the $0$-dimensional sphere
$S^0$, combinatorially is a subdivision of $\susp_1(v,K)$, we have that both spaces are PL-homeomorphic. In particular,
one-point suspensions provide an economic way of suspending a simplicial complex; see~\cite{BjoernerLutz2000,BjoernerLutz2003,Lutz2004a,Lutz2004b,Lutz2003gpre}.
\begin{figure}
\begin{center}
  \includegraphics[width=.7\linewidth]{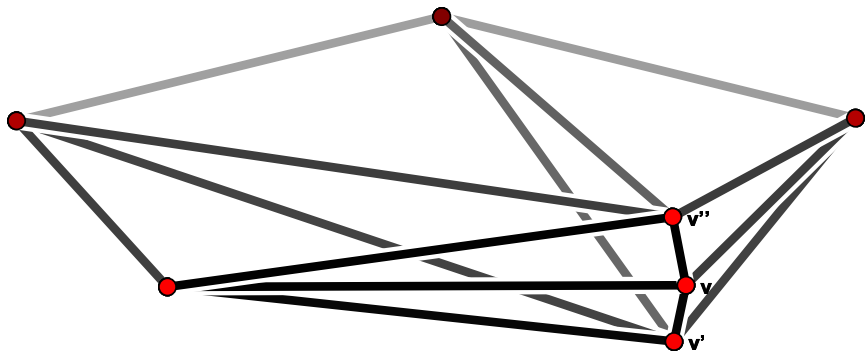}
\end{center}
\caption{The suspension of the circle $C_5$.}
\label{fig:circle5_1}
\end{figure}

\begin{figure}
\begin{center}
  \includegraphics[width=.7\linewidth]{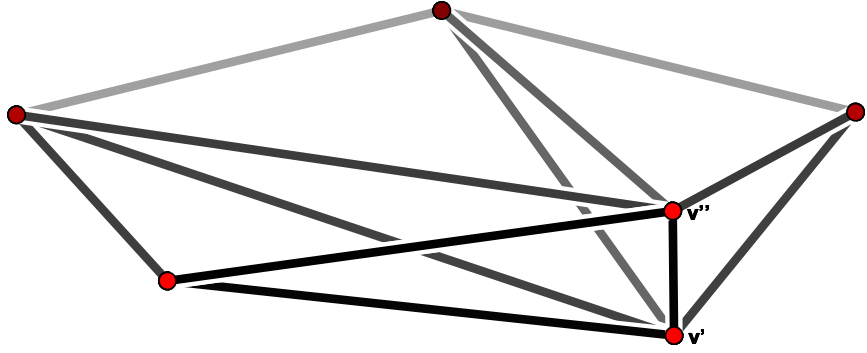}
\end{center}
\caption{The one-point suspension of the circle $C_5$ with respect to one of its vertices.}
\label{fig:circle5_2}
\end{figure}

\begin{exmp}
  Figure~\ref{fig:circle5_1} displays the $1$-skeleton of the suspension $S^0*C_5$ of the $5$-gon $C_5$. By ``removing''
  the original vertex $v$ from the suspension, we obtain the one-point suspension $\susp_1(v,C_5)$ of $C_5$ with respect to $v$;
  see Figure~\ref{fig:circle5_2}.
\end{exmp}

This one-point suspension has a higher-dimensional analog: Instead of the join product of a simplicial complex $K$ with
$S^0$, which is the boundary of an $1$-simplex, we can take the join product of $K$ with the boundary $\partial\Delta_d$
of a $d$-di\-men\-sional simplex $\Delta_d$ and then ``remove'' a vertex $v$ of $K$.
\begin{defn}
  Let $K$ be a simplicial complex and let $v$ be a vertex of $K$.  The \emph{reduced join} of $K$ with the boundary
  $\partial\Delta_d$ of a $d$-simplex $\Delta_d$ with respect to $v$ is the simplicial complex
  \[\partial\Delta_d*_v K=((\partial\Delta_d*K)\backslash (\partial\Delta_d*\star_K(v)))\cup (\Delta_d*\link_K(v)).\]
\end{defn}

From the construction of the reduced join we see that $\partial\Delta_d*_v K$ is obtained from $\partial\Delta_d*K$ by a
generalized bistellar flip which removes $\partial\Delta_d*\star_K(v)$ from $\partial\Delta_d*K$ and inserts
$\Delta_d*\link_K(v)$ instead (cf.~\cite{BjoernerLutz2000} for the definition of and further references on bistellar
flips).  The reverse direction of this operation is called \emph{starring a vertex in} $\Delta_d$ in \cite{BagchiDatta1998}.
Since $\partial\Delta_d*K$ is a subdivision of $\partial\Delta_d*_v K$, both spaces are PL-homeomorphic. In fact, the reduced
join $\partial\Delta_d*_v K$ can be described as $d$ iterated one-point suspensions of $K$ with respect to $v$ and copies of $v$
that are generated in each intermediate step.

The wreath product of polytopes has the following combinatorial analog.

\begin{defn}
  Let $K$ be a simplicial complex with $n$ vertices and let $\partial\Delta_d$ be 
  the $(d-1)$-dimensional boundary of an abstract $d$-simplex $\Delta_d$.
  We define the \emph{wreath product} $\partial\Delta_d\wr K$ of $\partial\Delta_d$ with $K$ as
  follows.  As vertices of $\partial\Delta_d\wr K$ we take $d+1$ copies $v_1^1,\dots,v_1^{d+1},\dots,v_n^1,\dots,v_n^{d+1}$ of
  the vertices $v_1,\dots,v_n$ of $K$. The facets of $\partial\Delta_d\wr K$ are all those subsets $S$ of vertices of
  $\partial\Delta_d\wr K$ of the form
  \[S:=\bigcup_{v\in F}\{v^1,\dots,v^{d+1}\}\cup\bigcup_{v\notin F}\{v^1,\dots,/,\dots,v^{d+1}\},\]
  where $F$ is a facet of $K$ and for the vertices $v\notin F$ 
  exactly one of the vertices $\{v^1,\dots,v^{d+1}\}$ is omitted.
\end{defn}

\begin{rem}\label{rem:wreath-facets}
  It follows from the construction that for $d>0$ every facet $S$ of $\partial\Delta_d\wr K$ arises from some facet $F$
  of $K$ as the multiple join product of copies of the full $d$-simplex $\Delta_d$ for every $v\in F$ with copies of
  facets of $\Delta_d$ for every $v\notin F$.  If $d=0$, then $\Delta_d$ is a point and $\partial\Delta_d\wr K=K$.  Also
  $\partial\Delta_d\wr K=\Delta_d$ if $K$ is a point.
\end{rem}

\begin{exmp}
  In Figure~\ref{fig:octa} we display one facet of the $14$-dimensional simplicial complex\,
  $\partial\Delta_2\wr\,\partial\,$octahedron\, that arises from the upper front triangle of the octahedron.  Every
  vertex of the upper front triangle contributes a full simplex $\Delta_2$ to the facet of\,
  $\partial\Delta_2\wr\,\partial\,$octahedron, all the other vertices contribute a $1$-dimensional maximal face of
  $\Delta_2$.
  \begin{figure}
    \begin{center}
      \includegraphics[width=.475\linewidth]{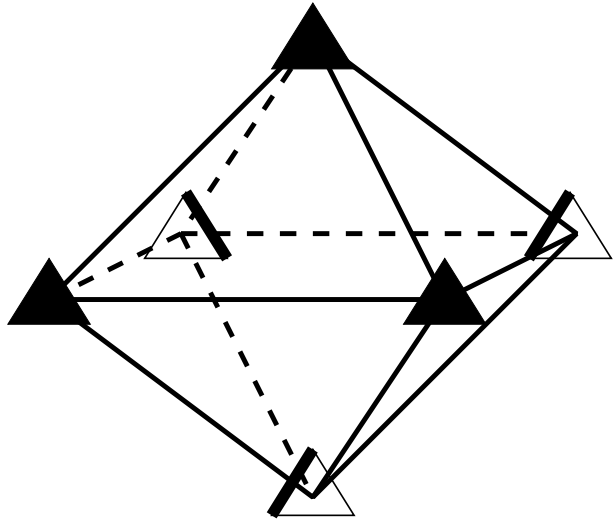}
    \end{center}
    \caption{One facet of\, $\partial\Delta_2\wr\,\partial\,$octahedron.}
    \label{fig:octa}
  \end{figure}
\end{exmp}

\begin{prop}\label{prop:reduced}
  Let $K$ be a simplicial complex with at least two (distinct) vertices $v_1,v_2$.  The reduced join is a commutative
  operation, i.e.,
  \[\partial\Delta_{d_1}*_{v_1}(\partial\Delta_{d_2}*_{v_2} K)\cong\partial\Delta_{d_2}*_{v_2}(\partial\Delta_{d_1}*_{v_1} K),\]
  for $d_1,d_2\ge0$.  In particular, $\partial\Delta_d\wr K$ can be obtained from $K$ by successive reduced joins (in an
  arbitrary order) with $\partial\Delta_d$ with respect to all the vertices of $K$.
\end{prop}
\begin{proof}
  Let $\{v_1,v_2,\dots,v_n\}$ be the set of vertices of~$K$.  
  As the vertices of $\partial\Delta_{d_1}*_{v_1}K$ we take $v_1^{1},\dots,v_1^{d_1+1},v_2,\dots,v_n$.  Then we have as
  facets of $\partial\Delta_{d_1}*_{v_1}K$, for all facets $F$ of $K$, all those subsets $S$ of vertices of
  $\partial\Delta_{d_1}*_{v_i}K$ of the form
  \[S=(F\backslash\{ v_1\})\cup\bigcup_{v_1\in F}\{ v_1^1,\dots,v_1^{d_1+1}\}
  \cup\bigcup_{v_1\notin F}\{ v_1^1,\dots,/,\dots,v_1^{d_1+1}\},\] where, if $v_1$ is not in $F$, exactly one of the
  vertices $\{ v_1^1,\dots,v_1^{d_1+1}\}$ is omitted.  
  As vertices of $\partial\Delta_{d_1}*_{v_1}(\partial\Delta_{d_2}*_{v_2} K)$ we take
  \[v_1^{1},\dots,v_1^{d_1+1},v_2^{1},\dots,v_2^{d_2+1},v_3,\dots,v_n.\]
  The facets of $\partial\Delta_{d_1}*_{v_1}(\partial\Delta_{d_2}*_{v_2} K)$ then are, for all facets $F$ of $K$, all those subsets $S$
  of vertices of $\partial\Delta_{d_1}*_{v_1}(\partial\Delta_{d_2}*_{v_2} K)$ of the form
  \begin{eqnarray*}
    S&=&(F\backslash\{ v_1,v_2\})\cup\bigcup_{v_1\in F}\{ v_1^1,\dots,v_1^{d_1+1}\}
    \cup\bigcup_{v_1\notin F}\{ v_1^1,\dots,/,\dots,v_1^{d_1+1}\}\\
    & & \hphantom{(F\backslash\{ v_1,v_2\})}   \cup\bigcup_{v_2\in F}\{ v_2^1,\dots,v_2^{d_2+1}\}
    \cup\bigcup_{v_2\notin F}\{ v_2^1,\dots,/,\dots,v_2^{d_2+1}\},
  \end{eqnarray*}
  where, if $v_1$ respectively $v_2$ is not in $F$, exactly one of the vertices $\{v_1^1,\dots,v_1^{d_1+1}\}$
  respectively $\{v_2^1,\dots,v_2^{d_2+1}\}$ is omitted.  The roles of $v_1$ and $v_2$ can clearly be exchanged, and
  hence the result follows.
\end{proof}

Similar to Proposition~\ref{prop:symmetry_wreath_product_polytope} for the corresponding polytopal construction, the
wreath product allows us to construct highly symmetric simplicial complexes.

\begin{prop}
  The wreath product of the automorphism groups $\Aut\partial\Delta_d\wr\Aut K=(S_{d+1})^n\rtimes\Aut K$, with respect
  to the natural action of $\Aut K$ on the $n$~vertices of~$K$, acts as a group of automorphisms of $\partial\Delta_d\wr K$.  
  In particular, if $\Aut K$ acts transitively on the set of vertices of $K$, 
  then also $\partial\Delta_d\wr K$ admits a vertex-transitive group
  of automorphisms.
\end{prop}

\mathversion{bold}
\subsection{$f$-vectors of Wreath Products}
\mathversion{normal} An $(e-1)$-dimensional simplicial complex $K$ is called \emph{pure} if all its maximal faces are of
dimension $e-1$.  Clearly, since the one-point suspensions have this property, the wreath product $\partial\Delta_d\wr K$ is
pure if and only if $K$ is pure.  The wreath product $\partial\Delta_d\wr K$ can be built from $nd$ iterated one-point
suspensions.  Since each one-point suspension step increases the dimension by one, we have that $\dim\partial\Delta_d\wr K=nd+e-1$.
  
Recall, that the \emph{$f$-vector} of the $(e-1)$-dimensional simplicial complex $K$ is the sequence
\[f(K)=(f_0,f_1,\dots,f_{e-1}),\] where $f_i$ is the number of $i$-dimensional faces of~$K$, for $0\leq i\leq e-1$.

\begin{prop}
  Abbreviating $n=f_0(K)$, the $f$-vector of\, $\partial\Delta_d\wr K$ has components
\begin{eqnarray*}
\lefteqn{f_i(\partial\Delta_d\wr K) = \sum_{j=\max\{ 0,i+1-nd\} }^{\min\{ e,\lfloor\frac{i+1}{d+1}\rfloor\}}
\left( f_{j-1}(K)\cdot
\sum_{\footnotesize\begin{array}{l}
      u_1\cdot 1+u_2\cdot 2+\dots +u_{d}\cdot d\\
      =i+1-j(d+1)\,\,\,\,\,\mbox{with}\\
      u_k\in{\mathbb N},\,\,\, 1\leq k\leq d
      \end{array}}\left[\binom{n-j}{u_{d}}\binom{d+1}{d}^{u_{d}}\cdot\right.
\vphantom{\sum_{\footnotesize\begin{array}{l}
      u_1\cdot 1+u_2\cdot 2+\dots +u_{d}\cdot d\\
      =i+1-j(d+1)\,\,\,\,\,\mbox{with}\\
      u_k\in{\mathbb N},\,\,\, 1\leq k\leq d
      \end{array}}}\right.}\\
&&\left.\vphantom{\sum_{\footnotesize\begin{array}{l}
      u_1\cdot 1+u_2\cdot 2+\dots +u_{d}\cdot d\\
      =i+1-j(d+1)\,\,\,\,\,\mbox{with}\\
      u_k\in{\mathbb N},\,\,\, 1\leq k\leq d
      \end{array}}}
\left.\binom{n-j-u_{d}}{u_{d-1}}\binom{d+1}{d-1}^{u_{d-1}}\dots\binom{n-j-u_{d}-u_{d-1}-\dots -u_2}{u_1}\binom{d+1}{1}^{u_1}\right]\right)
\end{eqnarray*}
for\, $0\leq i\leq nd+e-1$. In particular
\[f_0(\partial\Delta_d\wr K)=n(d+1)\]
and
\[f_{nd+e-1}(\partial\Delta_d\wr K)=f_{e-1}(K)(d+1)^{n-e}.\]
\end{prop}
\begin{proof}
  By definition, the vertex set of the complex $\partial\Delta_d\wr K$ is formed of $d+1$ copies 
  of the vertices of $K$. Hence, $f_0(\partial\Delta_d\wr K)=n(d+1)$.
  
  The facets of $\partial\Delta_d\wr K$ of dimension $nd+e-1$ arise from facets of $K$ of dimension $e-1$. For every
  $(e-1)$-dimensional facet $F$ of $K$ the corresponding facets of $\partial\Delta_d\wr K$ are of the form
  $\bigcup_{w\in F}\{ w^1,\dots,w^{d+1}\}\cup\bigcup_{w\notin F}\{ w^1,\dots,/,\dots,w^{d+1}\}$.  Since $F$ has
  cardinality~$e$, we take (a copy of) the full simplex $\Delta_d$ for the $e$ vertices in $F$ and (a copy of) a facet
  of $\Delta_d$ for the remaining $n-e$ vertices of $K$.  The simplex $\Delta_d$ has $(d+1)$ facets, thus there are
  $(d+1)^{n-e}$ facets of $\partial\Delta_d\wr K$ that arise from the $(e-1)$-dimensional facet $F$ of~$K$.  Moreover,
  as in Corollary~\ref{cor:number-of-facets}, $f_{nd+e-1}(\partial\Delta_d\wr K)=f_{e-1}(K)(d+1)^{n-e}$.
  
  Let $G$ be an $i$-dimensional face of $\partial\Delta_d\wr K$. Every vertex of $K$ either contributes a (copy of) a
  full simplex $\Delta_d$ or a (copy of) a face of $\Delta_d$ to $G$.  The set of vertices of $K$ that contribute a full
  simplex form a face $H$ of cardinality $j$ of $K$.  This face $H$ therefore contributes $j(d+1)$ vertices to $G$.
  Since $G$ is $i$-dimensional, there are $i+1-j(d+1)$ vertices of $G$ left that are contributed to by the other $n-j$
  vertices of $K$.  In fact, every of the $n-j$ vertices contributes between $1$ and $d$ vertices to
  $\partial\Delta_d\wr K$, so let $u_k$ be the number of vertices of $K$ that contribute $k$ vertices to $G$.  Since $G$
  is $i$-dimensional, it follows that $u_1\cdot 1+u_2\cdot 2+\dots +u_{d}\cdot d=i+1-j(d+1)$.  There are
  $\binom{n-j}{u_{d}}\binom{d+1}{d}^{u_{d}}$ choices for $u_{d}$ of the $n-j$ vertices to contribute $d$ vertices to
  $\partial\Delta_d\wr K$, etc. Altogether, there are $f_{j-1}(K)$ faces $H$ of cardinality $j$ of $K$ that can
  contribute for each vertex a full simplex $\Delta_d$ to an $i$-dimensional face $G$ of $\partial\Delta_d\wr K$.
  Observe, that $j$ has to be restricted to the range\, $\max\{ 0,i+1-nd\} \leq j\leq\min\{
  e,\lfloor\frac{i+1}{d+1}\rfloor\}$.
\end{proof}

A simplicial complex $K$ is called \emph{$k$-neighborly} if $f_i(K)=\binom{f_0(K)}{i+1}$ 
for $0\leq i\leq k-1$, that is, every set of $k$ (or less) vertices is a face of $K$.

\begin{prop} 
\label{prop:neighborly}
If $K$ is a $k$-neighborly simplicial complex, then $\partial\Delta_d\wr K$
is ($k(d+1)+d$)-neighborly. If $K$ is not $k$-neighborly, then $\partial\Delta_d\wr K$
is not $k(d+1)$-neighborly.
\end{prop}
\begin{proof}
Let $K$ be $k$-neighborly and let $F$ be a set of vertices of 
$\partial\Delta_d\wr K$ of cardinality $k(d+1)+d$.
Every vertex of $K$ contributes at least $d$ and at most $d+1$ vertices
to every facet of $\partial\Delta_d\wr K$. Since $F$ has cardinality $k(d+1)+d$,
there are at most $k$ vertices of $K$ for which all its $d+1$ copies
are present in $F$. However, $K$ is $k$-neighborly, so there is indeed
a facet of $K$ that contains these at most $k$ vertices. One of the corresponding
facets of $\partial\Delta_d\wr K$ then contains $F$.

Let $K$ be not $k$-neighborly, and suppose that every set of vertices of 
$\partial\Delta_d\wr K$ of cardinality $k(d+1)$ is a face of $\partial\Delta_d\wr K$.
Since $K$ is not $k$-neighborly, there is a set
$G$ of $k$ vertices of $K$ that is not a face of $K$. The union of 
the $d+1$ copies of these $k$ vertices then is a set of cardinality $k(d+1)$
which is not a face of $\partial\Delta_d\wr K$. Contradiction.
\end{proof}

\section{Combinatorial Decompositions of One-point Suspensions and Wreath Products}
\emph{Vertex-decomposability}, \emph{shellability}, and \emph{constructibility} are three standard concepts to decompose
a pure simplicial complex into its collection of facets; see Bj\"orner\cite{Bjoerner1995}.  We show that these properties
are respected by one-point suspensions and hence also by the wreath product construction.  A pure $(e-1)$-dimensional
simplicial complex $K$ is

\begin{itemize}
\item \emph{vertex-decomposable} if either $K$ is a simplex (possibly $\{\emptyset\}$) or there is a vertex~$v$ such
  that the link $\link_{K}(v)$ and the deletion $\del_{K}(v)$ of $v$ in $K$ are both vertex-decomposable simplicial
  complexes;
  
\item \emph{shellable} if it has a \emph{shelling}, i.e., there is a linear ordering $F_1,F_2,\dots,F_{f_{e-1}(K)}$ of
  the $f_{e-1}(K)$ facets of $K$ such that $(2^{F_1}\cup\dots\cup 2^{F_{k-1}})\cap 2^{F_k}$ is a pure
  $(e-2)$-dimensional simplicial complex for $2\leq k\leq f_{e-1}(K)$, where $2^{F}$ is the set of all faces of a
  simplex~$F$;
  
\item \emph{constructible} if either $K$ is a simplex or there are two $(e-1)$-dimensional constructible subcomplexes
  $K_1$ and $K_2$ of $K$ such that their union is $K$ and their intersection is an $(e-2)$-dimensional constructible
  simplicial complex;

\item \emph{Cohen-Macaulay} (with respect to some field $\FF$) if the reduced homology groups
  $\tilde{H}_i(\link_{K}(G);\FF)$ vanish for $i\neq \dim(\link_{K}(G))$ for all faces $G\in K$.
\end{itemize}
For pure simplicial complexes the following implications are strict (cf.~\cite{Bjoerner1995}):
\[\text{vertex-decomposable} \,\Rightarrow\, \text{shellable} \,\Rightarrow\, \text{constructible}
\,\Rightarrow\, \text{Cohen-Macaulay}.\]

Note that, due to Munkres~\cite{Munkres1984b}, Cohen-Macaulayness over a field is not a combinatorial property but an
entirely topological one.  Here we mention it for systematic reasons.  Munkres result~\cite{Munkres1984b} already implies
that the one-point suspension of a Cohen-Macaulay complex (and hence also any wreath product) is again Cohen-Macaulay.
Conversely, Cohen-Macaulayness of~$K$ is necessary for the Cohen-Macaulayness of~$\susp_1(v,K)$ since $K$ occurs as a link.

\begin{prop}{\rm (Provan and Billera~\cite[Proposition 2.5]{ProvanBillera1980})}
  The one-point suspension $\susp_1(v,K)$ is vertex-decomposable if and only if $K$ is.
\end{prop}

\begin{cor}\label{cor:vertex-decomposable}
  The wreath product $\partial\Delta_d\wr K$ is vertex-decomposable if and only $K$ is.
\end{cor}

\begin{prop}\label{prop:shellable}
  The one-point suspension $\susp_1(v,K)$ is shellable if and only if $K$ is.
\end{prop}
\begin{proof}
  Let $F_1,\dots,F_{f_{e-1}(K)}$ be a shelling order of the facets of~$K$.  As pointed out in Remark~\ref{rem:susp-facets} we have
  three kinds of facets in $\susp_1(v,K)$.  Thus we obtain a shelling order of the facets of $\susp_1(v,K)$ by replacing
  each facet $F_i$ which does not contain~$v$ by the pair $v'*F_i,v''*F_i$ and each facet $F_j$ which contains~$v$ by the
  facet $\overline{v'v''}*(F_j\setminus\{v\})$.

  For the converse observe that under the map $\beta$, defined in Remark~\ref{rem:susp-facets}, each shelling order of
  the facets of~$\susp_1(v,K)$ also induces a shelling order of the facets of~$K$ (after removing doubles).
\end{proof}

\begin{exmp}
  Iteratively applying the construction in the proof of Proposition~\ref{prop:shellable} yields shellings of wreath
  products $\partial\Delta\wr K$ from shellings of~$K$.
  
  In particular, if $F_1,F_2,\dots,F_{f_{e-1}(K)}$ is a shelling of $K$, then we first partition the facets of
  $\partial\Delta_d\wr K$ into $f_{e-1}(K)$ sets of facets $B(F_k)$ that arise from the facets $F_k$, $1\leq k\leq
  f_{e-1}(K)$, according to Remark~\ref{rem:wreath-facets}.  Each collection $B(F_k)$ is a join product $*_{w\in
    F_k}\Delta_d*_{w\notin F_k}\partial \Delta_d$, and therefore it is a shellable ball. For an explicit shelling of the
  first ball $B(F_1)$ we start with some of its facets and continue with those facets in $B(F_1)$ that differ from the
  first facet by two vertices, then with those facets that differ by four vertices, etc.
  
  As an example, we display in Figure~\ref{fig:shelling} a corresponding shelling of the set of facets associated with
  the upper front triangle of\, $\partial\Delta_2\wr\,\partial\,$octahedron\, from Figure~\ref{fig:octa}.

  \begin{figure}[htbp]
    \begin{center}
      \includegraphics[width=.35\linewidth]{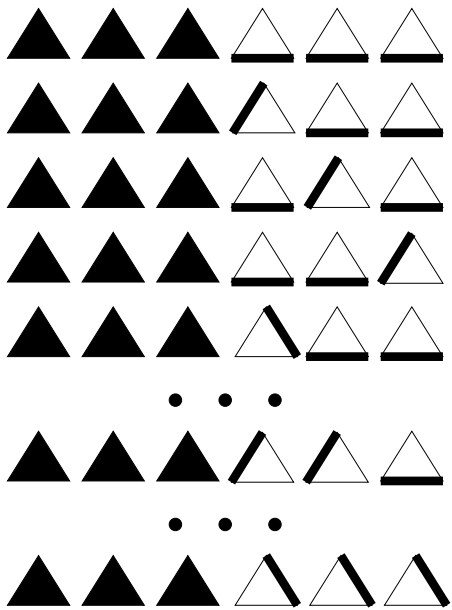}
    \end{center}
    \caption{A shelling of the facets associated with the upper front
      triangle of\, $\partial\Delta_2\wr\,\partial\,$octahedron.}
    \label{fig:shelling}
  \end{figure}
  
  The way that we have chosen the facets, we ensure that for every new facet in the ordering the intersection with the
  previous facets is\, ($nd+e-2$)-dimensional.  Upon completion of the shelling of $B(F_1)$ we continue with the facets
  of $B(F_2)$, etc; see Figure~\ref{fig:shelling2}.

  \begin{figure}[htbp]
    \begin{center}
      \includegraphics[width=.35\linewidth]{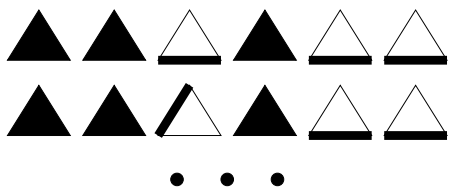}
    \end{center}
    \caption{The shelling of a consecutive collection of facets of\,
      $\partial\Delta_2\wr\,\partial\,$octahedron.}
    \label{fig:shelling2}
  \end{figure}
\end{exmp}

\begin{cor}\label{cor:shellable}
  The wreath product $\partial\Delta_d\wr K$ is shellable if and only if $K$ is shellable.
\end{cor}

\begin{prop}
  The one-point suspension $\susp_1(v,K)$ is constructible if and only if $K$ is.
\end{prop}
\begin{proof}
  A construction order of a simplicial complex is a sequence of increasingly fine (special) equivalence relations on the
  set of facets such that the final equivalence relation is the identity.  Clearly, for each facet $F$ of $K$ (not
  containing~$v$) the simplicial complex of two facets $v'*F$ and $v''*F$ is a constructible ball.  So, by virtue of the
  inverse $\beta^{-1}$ of the blocking map, which maps sets of facets of~$K$ to sets of facets of $\susp_1(v,K)$, and by
  an obvious induction on the dimension of~$K$, each construction of $K$ induces a construction of the one-point
  suspension.
  
  Conversely, if a pure simplicial complex is constructible, then all its vertex-links are constructible (see \cite{Bjoerner1995}
  and \cite{HachimoriZiegler2000}). Since $K$ appears as a vertex-link in every one-point suspension of $K$, a one-point
  suspension of $K$ is non-constructible if $K$ is non-constructible.
\end{proof}

\begin{cor}
\label{cor:constructible}
  The wreath product $\partial\Delta_d\wr K$ is constructible if and only if $K$ is constructible.
\end{cor}

\subsection{Combinatorial Strengthenings and Topological Weakenings of Contractibility}

Combinatorial notions which imply contractibility appear in various contexts in topology and combinatorics.  An
$(e-1)$-dim\-en\-sion\-al simplicial complex $K$ is
\begin{itemize}
\item \emph{non-evasive} if either $K$ is a single point or there is a vertex $v$ of $K$ such that both $\link_{K}(v)$
  and $\del_{K}(v)$ of $v$ are non-evasive;
\item \emph{collapsible} if the Hasse diagram of~$K$ (seen as a graph whose edges are directed towards the
  higher-dimensional faces, and $\emptyset$ counts as a face of~$K$) admits a perfect matching which is \emph{acyclic},
  that is, the graph remains acyclic if the orientations of the edges in the matching are reversed;
\item \emph{contractible} if $K$ is homotopy equivalent to a point;
\item \emph{$\ZZ$-acyclic} if all reduced homology groups of~$K$ with integer coefficients vanish.
\end{itemize}

For simplicial complexes the following implications are strict (cf.\ \cite{Bjoerner1995}):
\[\mbox{\rm cone} \,\Rightarrow\, \mbox{\rm non-evasive} \,\Rightarrow\, \mbox{\rm collapsible}
\,\Rightarrow\, \mbox{\rm contractible}\,\Rightarrow\, \mbox{\rm $\ZZ$-acyclic} \,\Rightarrow\, \mbox{\rm
  $\tilde{\chi}=0$} ,\] where $\tilde{\chi}$ denotes the reduced Euler characteristics of a simplicial complex.  The
perfect matching in the definition of collapsibility is a special case of a Morse matching in the sense of
Chari~\cite{Chari2000}; see also Forman~\cite{Forman1998,Forman2002}.  In the sequel we call a perfect acyclic
matching a \emph{perfect Morse matching} and the unique vertex matched to the empty face is called \emph{critical}.  The
concept of evasiveness originally stems from the complexity theory of graph properties and was reformulated in terms of
simplicial complexes by Kahn, Saks, and Sturtevant \cite{KahnSaksSturtevant1984}; see also \cite{Bjoerner1995},
\cite{Lutz2002}, and \cite{Welker1999}.

\begin{prop}
  The one-point suspension $\susp_1(v,K)$ is a cone if and only if $K$ is.
\end{prop}
\begin{proof}
  Suppose that $K=a*B$ is a cone with apex~$a$.  To prove that $\susp_1(v,a*B)$ is a cone we distinguish two cases: If
  $v=a$, then $\susp_1(a,a*B)=a'*(a''*B)$.  If $v\ne a$, then $v\in B$ and
  $\susp_1(v,a*B)=(a*\partial\overline{v'v''}*B\setminus\partial\overline{v'v''}*a*\star_B(v))\cup\overline{v'v''}*a*\link_B(v)=a*\susp_1(v,B)$.
  
  For the converse assume that $\susp_1(v,K)=a*B$ for some vertex $a$ and some induced subcomplex~$B$.  If $a=v'$ (or,
  symmetrically, $a=v''$) then $\susp_1(v,K)=\overline{v'v''}*C$ where $C$ is the subcomplex of $\susp_1(v,K)$ induced
  on the complement of $\{v',v''\}$.  Clearly, $C=\link_K(v)$ and $K=v*C$ is a cone.  Otherwise if $a\not\in\{v',v''\}$
  then $a\in K\setminus\star_K(v)$ and $K=\star_K(a)\cup\link_K(a)$ is a cone with apex~$a$.
\end{proof}

\begin{cor}
  The wreath product $\partial\Delta_d\wr K$ a cone if and only if $K$ is.
\end{cor}

\begin{prop}
  The one-point suspension $\susp_1(v,K)$ is non-evasive if and only if $K$ is.
\end{prop}
\begin{proof}
  Let $K$ be non-evasive. Then clearly, $\link_{\susp_1(v,K)}(v')=K$ is non-evasive by assumption.  Furthermore,
  $\del_{\susp_1(v,K)}(v')$ is a cone with apex $v''$, and therefore it is also non-evasive.
  
  If $\susp_1(v,K)$ is non-evasive, then there is a vertex $w$ such that $\link_{\susp_1(v,K)}(w)$ and
  $\del_{\susp_1(v,K)}(w)$ are non-evasive.  If $w\in\{ v',v''\}$, then it follows that $\link_{\susp_1(v,K)}(w)=K$ is
  non-evasive.  Thus, let us assume that $w\notin\{ v',v''\}$. In this case,
  $\link_{\susp_1(v,K)}(w)=\susp_1(v,\link_K(w))$ and $\del_{\susp_1(v,K)}(w)=\susp_1(v,\del_K(w))$ are non-evasive, so
  by induction, $K$ is non-evasive.
\end{proof}

\begin{cor}
  The wreath product $\partial\Delta_d\wr K$ is non-evasive if and only if $K$ is.
\end{cor}

\begin{prop}\label{prop:collapse}
  If $K$ is collapsible, then the one-point suspension $\susp_1(v,K)$ is collapsible.
\end{prop}
\begin{proof}
  We prove that each perfect Morse matching $\mu$ of $K$ can be lifted to a perfect Morse matching $\bar\mu$ of
  $\susp_1(v,K)$.  This lifting is not canonical but it depends on choices.
  
  Let $(\sigma,\tau)\in\mu$.  Depending on the relative positions of the faces $\sigma,\tau$ to the special vertex~$v$
  they may induce up to three different matched pairs in $\bar\mu$, as it will be defined now.  We distinguish the
  following cases:
  \begin{enumerate}
  \item $\sigma,\tau\in\star_K(v)$: Then we let
    $(\sigma\setminus\{v\}\cup\{v',v''\},\tau\setminus\{v\}\cup\{v',v''\})\in\bar\mu$.\label{mu:1}
  \item $\sigma\in\star_K(v)$ and $\tau\not\in\star_K(v)$: In this case we necessarily have $\sigma=\tau\cup\{v\}$, and
    we let $(\sigma\setminus\{v\}\cup\{v'\},\tau),\ 
    (\sigma\setminus\{v\}\cup\{v',v''\},\tau\cup\{v''\})\in\bar\mu$.\label{mu:2}
  \item $\sigma,\tau\not\in\star_K(v)$: Then we let
    $(\sigma,\tau),\ (\sigma\cup\{v'\},\tau\cup\{v'\}),\ (\sigma\cup\{v''\},\tau\cup\{v''\})\in\bar\mu$.\label{mu:3}
  \end{enumerate}
  
  In order to prove that $\bar\mu$ is indeed a perfect matching in the Hasse diagram of~$\susp_1(v,K)$ we cannot avoid a
  somewhat tedious case distinction according to the six different types of pairs in~$\bar\mu$ which we address as
  \ref{mu:1}, \ref{mu:2}a, \ref{mu:2}b, \ref{mu:3}a, \ref{mu:3}b, and \ref{mu:3}c, respectively: Let
  $\phi\in\susp_1(v,K)$ be any face.
  
  \framebox{$\phi\in K\setminus\star_K(v)$:} Let $\psi$ be the match of~$\phi$ in~$\mu$.  If $\psi\in
  K\setminus\star_K(v)$, too, then (type~\ref{mu:3}a) $\psi$ is the unique match of~$\phi$ in~$\bar\mu$.  Otherwise
  $\psi=\phi\cup\{v\}$ and $(\phi,\psi\setminus\{v\}\cup\{v'\})\in\bar\mu$ (type~\ref{mu:2}a).
  
  \framebox{$\phi=\phi'\cup\{v'\}$ and $\phi'\in K\setminus\link_K(v)$:} Let $(\phi',\psi')\in\mu$.  Then
  $\psi'\not\in\star_K(v)$, and the unique match of~$\phi$ is $\psi'\cup\{v'\}$ (type~\ref{mu:3}b).

  \framebox{$\phi=\phi'\cup\{v'\}$ and $\phi'\in\link_K(v)$:} Again let $(\phi',\psi')\in\mu$.  If
  $\psi'\not\in\star_K(v)$, then (type~\ref{mu:3}b) $\psi'\cup\{v'\}$ is the unique match.  Otherwise $\psi'\in\star_K(v)$
  and hence $\phi'\cup\{v\}=\psi'$, and thus (type~\ref{mu:2}a) the unique match is $\phi'$.

  \framebox{$\phi=\phi''\cup\{v''\}$ and $\phi''\in K\setminus\link_K(v)$:} Let $(\phi'',\psi'')\in\mu$.  Then
  $\psi''\not\in\star_K(v)$ and (type~\ref{mu:3}c) the unique match of~$\phi$ is $\psi''\cup\{v''\}$.

  \framebox{$\phi=\phi''\cup\{v''\}$ and $\phi''\in\link_K(v)$:} Again let $(\phi'',\psi'')\in\mu$.  If
  $\psi''\not\in\star_K(v)$, then (type~\ref{mu:3}c) $\psi''\cup\{v''\}$ is the unique match of $\phi$ in~$\bar\mu$.
  Otherwise $\psi''\cup\{v\}=\phi''$, and the unique match is $\phi''\setminus\{v\}\cup\{v',v''\}$ (type \ref{mu:2}b).
  
  \framebox{$\phi=\phi'''\cup\{v',v''\}$:} Then $\phi'''\in\link_K(v)$.  Let $(\phi'''\cup\{v\},\psi''')\in\mu$.  If
  $\psi'''\in\star_K(v)$, then $\phi$ is matched to $\psi'''\setminus\{v\}\cup\{v',v''\}$ (type \ref{mu:1}).  Otherwise
  $\psi'''\in\link_K(v)$ and hence $\psi'''=\phi'''$.  This is the case \ref{mu:2}b, and we conclude that
  $(\phi,\phi'''\cup\{v''\})\in\bar\mu$.

  The acyclicity of $\bar\mu$ is inherited from the acyclicity of~$\mu$; we omit the details.
\end{proof}

It seems to be an open question whether the converse of the previous proposition holds.  However, there exist perfect
Morse matchings of one-point suspensions which are not induced by perfect Morse matchings of the base space.

\begin{exmp}
  Let $\pi$ be the $1$-dimensional simplicial complex on the vertex set $\{1,2,3,4\}$ with facets $12,23,34$, that is,
  $\pi$ is a path on four vertices. In the Hasse diagram of~$\pi$ we consider the perfect Morse matching
  $\mu=\{(12,1),(23,2),(34,3),(4,\emptyset)\}$.  Figure~\ref{fig:collapse} displays $\mu$ and its lifting $\bar\mu$ as
  defined in the proof of Proposition~\ref{prop:collapse}.
  \begin{figure}[htbp]\centering
    \subfigure[The perfect Morse matching $\mu$ of the path~$\pi$.  The vertex~$4$ is
    critical.\label{fig:mu}]{\includegraphics[width=.45\textwidth]{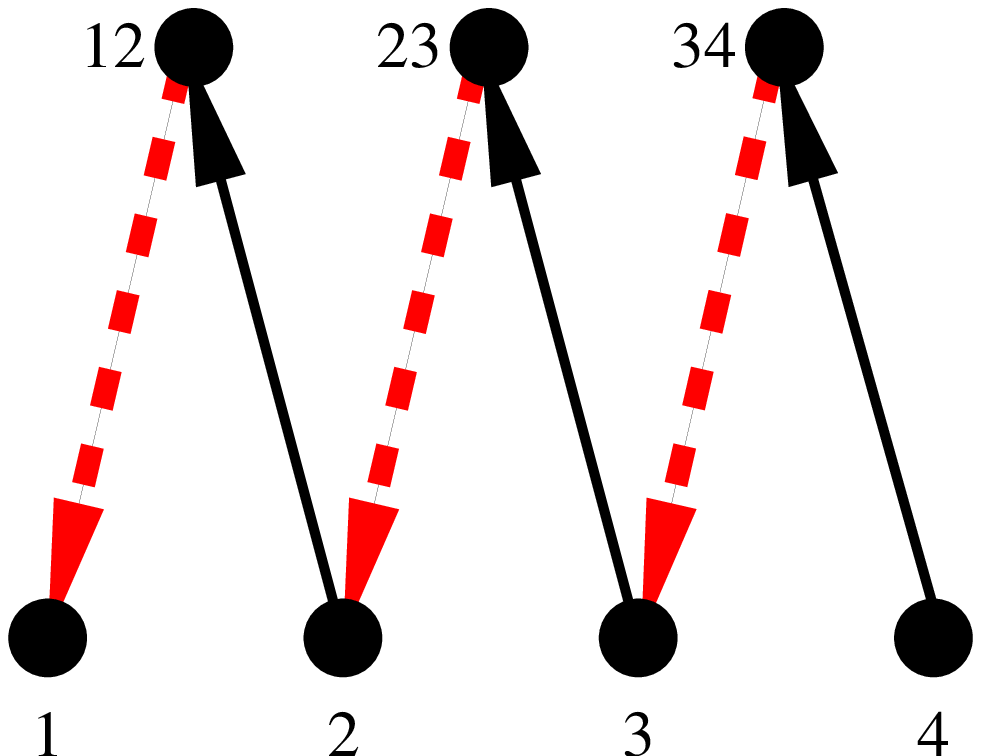}}\quad \subfigure[The perfect Morse
    matching $\bar\mu$ of $\susp(4,\pi)$.  The vertex~$4'$ is
    critical.\label{fig:bar-mu}]{\includegraphics[width=.45\textwidth]{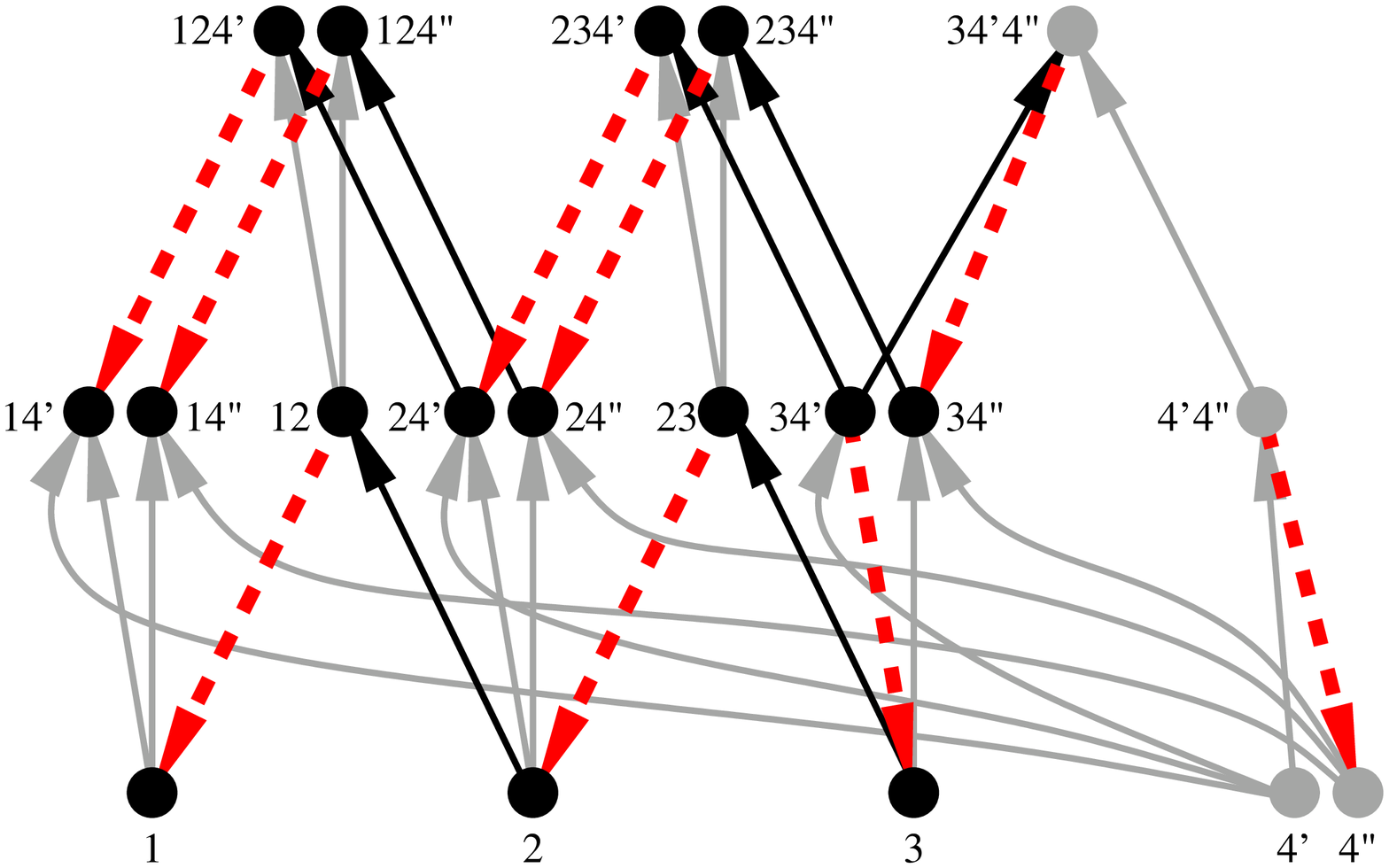}}
    \caption{The dashed arrows pointing downwards form the respective matchings.  In both cases the empty face is
      omitted.\label{fig:collapse}}
  \end{figure}
\end{exmp}

\begin{cor}
  If $K$ is collapsible, then the wreath product $\partial\Delta_d\wr K$ is collapsible.
\end{cor}

The homology of a suspension is the same as the (reduced) homology of the base space, 
up to a shift in dimension: $\widetilde H_i(\susp_1(v,K))=\widetilde H_{i-1}(K)$, for $i\ge1$.
Hence, the one-point suspension $\susp_1(v,K)$ is $\ZZ$-acyclic if and only 
if $K$ is $\ZZ$-acyclic, and likewise for the wreath products.
Since the suspension of a $\ZZ$-acyclic space is even contractible (cf.\ \cite{Bjoerner1995} and \cite{Lutz2002}) 
we have the following stronger result.
\begin{prop}
If $K$ is $\ZZ$-acyclic, then the one-point suspension $\susp_1(v,K)$ is contractible.
\end{prop}
\begin{cor}
  If $K$ is $\ZZ$-acyclic, then the wreath product $\partial\Delta_d\wr K$ is contractible for $d\geq 1$.
\end{cor}

As vertex-transitivity translates to wreath products, the wreath product $\partial\Delta_d\wr K$
of a vertex-transitive $\ZZ$-acyclic simplicial complex $K$ yields for $d\geq 1$
a vertex-transitive contractible simplicial complex.
A first example of a vertex-transitive $\ZZ$-acyclic simplicial
complex was constructed by Oliver; for further examples, based on the
$2$-skeleton of the Poincar\'e homology $3$-sphere in its description
by Threlfall and Seifert \cite{ThrelfallSeifert1931} and 
Weber and Seifert \cite{WeberSeifert1933} as the 
\emph{spherical dodecahedron space}, see \cite{Lutz2002}.  
In particular, the smallest currently known $\ZZ$-acyclic vertex-transitive simplicial complex
is the $5$-dimensional complex $K_3$ with $30$ vertices of Lutz \cite{Lutz2002}.
\begin{thm}
The wreath products $\partial\Delta_d\wr K_3$ of the vertex-transitive $5$-dimensional $\ZZ$-acyclic simplicial
complex $K_3$ with $30$ vertices give $(30d+5)$-dimensional vertex-transitive contractible complexes 
with $30(d+1)$ vertices for $d\geq 1$.
\end{thm}
\begin{rem}
The previously known vertex-transitive contractible simplicial complexes with $60$ vertices
from \cite{Lutz2002} are of dimension $11$, $23$, and $29$. Therefore, the wreath product
construction provides new examples. These are of particular interest, since it is still
open whether there are vertex-transitive collapsible simplicial complexes, and, if such spaces exist,
whether they can be constructed by starting with contractible or $\ZZ$-acyclic vertex-transitive complexes.
Non-existence would, on the other hand, prove the long-standing evasiveness conjecture 
for graph properties; cf.\ \cite{KahnSaksSturtevant1984}.
\end{rem}

\section{PL-Topology of Wreath Products}

A pure $(e-1)$-dimensional simplicial complex $K$ is a \emph{(weak) simplicial pseudomanifold (with boundary)} if every
$(e-2)$-dimensional face is contained in exactly two (at most two) $(e-1)$-dimensional facets of the complex $K$.  Since
the one-point suspension of a space is PL-equivalent to the ordinary suspension it is clear that the one-point
suspension is a pseudomanifold if and only if the base space is.  Again this property extends to wreath products.

\subsection{Wreath Products of Spheres}

As wreath products are iterated one-point suspensions it is clear that the wreath product $\partial\Delta_d\wr S$ of a simplicial
sphere~$S$ is again a simplicial sphere.  This section is devoted to the study of how additional structures on~$S$
behave with respect to wreath products.

Newman \cite{Newman1926c} proved that a constructible pseudomanifold (with boundary) is a 
PL sphere (PL ball); see also Bj\"orner \cite{Bjoerner1995}.
In Propositions~\ref{cor:vertex-decomposable}, \ref{cor:shellable}, and \ref{cor:constructible}, we already proved that
wreath products $\partial\Delta_d\wr K$ inherit vertex-decomposability, shellability, or constructibility from the corresponding
property of the base space~$K$.

For simplicial pseudomanifolds we have the following implications: 

\begin{center}
\small
\noindent
\begin{picture}(160,160)(0,20)
\put (0,170) {polytopal sphere}
\put (42,158) {$\Searrow$}
\put (116,158) {$\Swarrow$}
\put (103,170) {vertex-decomposable sphere}
\put (55,146) {shellable sphere}
\put (80,134) {$\Downarrow$}
\put (49,122) {constructible sphere}
\put (80,110) {$\Downarrow$}
\put (47,98) {combinatorial sphere}
\put (80,86) {$\Downarrow$}
\put (53,74) {simplicial sphere}
\put (80,62) {$\Downarrow$}
\put (52.5,50) {homology sphere.}
\put (80,38) {$\Downarrow$}
\put (38,26) {Cohen-Macaulay complex}
\end{picture}
\end{center}

\begin{prop}
  The one-point suspension $\susp_1(v,S)$ is a polytopal sphere if and only if $S$ is.
\end{prop}
\begin{proof}
  If $S=\partial P$ for some simplicial polytope~$P$, then $\susp_1(v,S)\cong\partial\dw(v,P)$ is polytopal, too; see
  Section~\ref{sec:dual-wedge}.  In order to prove the converse, suppose that $\susp_1(v,S)=\partial Q$ for some
  simplicial polytope~$Q$.  Then the vertex figure $P/v'$ of the vertex $v'$ is a simplicial polytope whose boundary is
  isomorphic to~$S$ (as a simplicial complex).
\end{proof}

\begin{cor}\label{cor:polytopal}
  The wreath product $\partial\Delta_d\wr S$ is a polytopal sphere if and only if $S$ is.
\end{cor}

A simplicial $(e-1)$-sphere $K$ is a \emph{combinatorial} sphere if $K$ is PL-homeomorphic to the boundary of the
standard $e$-simplex $\Delta_e$.  In particular, all vertex-links of a combinatorial $(e-1)$-sphere are combinatorial
$(e-2)$-spheres. Observe that in all dimensions $e-1\neq 4$, every $(e-1)$-simplicial sphere with the property that all
its vertex-links are combinatorial $(e-2)$-spheres is itself a combinatorial sphere. In dimension $e-1=4$, it is an open
problem whether exotic simplicial $4$-spheres exist that are not combinatorial, but for which all vertex-links are
combinatorial $3$-spheres.  Since one-point suspensions are PL-equivalent to ordinary suspensions, one-point suspensions
(and thus also wreath products) of combinatorial spheres are again combinatorial spheres, and conversely.  

Vertex-transitive triangulations of combinatorial spheres with up to $15$ vertices (except for some
small symmetry groups) were enumerated in \cite{Lutz1999}. Among these, various examples are
polytopal wreath product spheres.
\begin{thm}\label{thm:enumeration_examples}
The vertex-transitive wreath-product spheres with $n(d+1)\leq 15$ are the following:
The spheres $\mbox{}^4\hspace{.3pt}6^{\,16}_{\,1}$, $\mbox{}^5\hspace{.3pt}8^{\,47}_{\,1}$,
$\mbox{}^6\hspace{.3pt}10^{\,23}_{\,1}$, $\mbox{}^7\hspace{.3pt}12^{\,193}_{\,1}$,
and $\mbox{}^8\hspace{.3pt}14^{\,38}_{\,1}$ of dimension $4$, $5$, $6$, $7$, and $8$
with $6$, $8$, $10$, $12$, and $14$ vertices from \cite{Lutz1999} are the polytopal wreath product
spheres $\partial\Delta_1\wr\partial C_2(n)$, $3\leq n\leq 7$, respectively.
The spheres $\mbox{}^7\hspace{.3pt}9^{\,34}_{\,1}$, $\mbox{}^9\hspace{.3pt}12^{\,299}_{\,1}$,
and $\mbox{}^{11}\hspace{.3pt}15^{\,86}_{\,1}$ of dimension $7$, $9$, and $11$ with
$9$, $12$, and $15$ vertices from \cite{Lutz1999} are the polytopal wreath product
spheres $\partial\Delta_2\wr\partial C_2(n)$, $3\leq n\leq 5$, respectively.
We have the identities
$\partial\Delta_k\wr\partial\Delta_2=\partial\Delta_2\wr\partial\Delta_k=\partial\Delta_{3k+2}$,
for $2\leq k\leq 4$, and
$\partial\Delta_k\wr\partial\Delta_1=\partial\Delta_1\wr\partial\Delta_k=\partial\Delta_{2k+1}$,
for $1\leq k\leq 6$.
Moreover, in the notation of \cite{Lutz1999}, 
$\mbox{}^8\hspace{.3pt}12^{\,294}_{\,1}=\partial\Delta_1\wr\partial C_3^{\Delta}=(\partial\Delta_3)^{*3}$,
where $C_3^{\Delta}$ is the $3$-dimensional cross-polytope;
and for the cyclic $4$-polytopes $C_4(6)$ and $C_4(7)$ we have that
$\mbox{}^9\hspace{.3pt}12^{\,299}_{\,1}=\partial\Delta_1\wr\partial C_4(6)=\partial C_{10}(12)=(\partial\Delta_5)^{*2}$\, and\,
$\mbox{}^{10}\hspace{.3pt}14^{\,38}_{\,1}=\partial\Delta_1\wr\partial C_4(7)$.
\end{thm}
\begin{rem}
Mani \cite{Mani1972} proved that every simplicial $(nd+e-1)$-sphere 
with $n(d+1)$ vertices is polytopal whenever $n(d+1)\leq (nd+e-1)+4$,
which settles the polytopality for most of the spheres of Theorem~\ref{thm:enumeration_examples}.
However, for the polytopality of the spheres $\mbox{}^7\hspace{.3pt}12^{\,193}_{\,1}$
and $\mbox{}^8\hspace{.3pt}14^{\,38}_{\,1}$ the characterization
as wreath products of polytopal spheres according to Corollary~\ref{cor:polytopal}
is needed.
\end{rem}

A (simplicial) \emph{homology $(e-1)$-sphere} is a manifold with the homology of the standard sphere
$S^{e-1}$.

\begin{prop}
  If $K$ is a homology $(e-1)$-sphere, different from the standard sphere $S^{e-1}$, and $d\geq 1$, 
  then $\partial\Delta_d\wr K$ is a non-PL sphere.
\end{prop}
\begin{proof}
  By the double suspension theorem of Edwards~\cite{Edwards1975} for the Mazur homology $3$-sphere and its generalization to
  arbitrary homology spheres by Cannon~\cite{Cannon1979}, the double suspension of every simplicial homology sphere $K$
  (different from the standard sphere) is a non-PL sphere. 
  If $d\geq 1$, then already the double reduced join $\partial\Delta_{d}*_{v_i}(\partial\Delta_{d}*_{v_j} K)$ with
  respect to two distinct vertices $v_i$ and $v_j$ of $K$ and therefore also the iterated reduced join $\partial\Delta_d\wr K$
  are PL homeomorphic to join products of spheres with the double suspension of $K$ and thus are simplicial spheres.
  They are non-PL spheres, since the homology sphere $K$ appears as the link of some of their faces.
\end{proof}

\begin{cor}
If $K$ is a vertex-transitive non-spherical homology sphere and $d\geq 1$,
then $\partial\Delta_d\wr K$ is a non-PL sphere for $d\geq 1$.
\end{cor}

Examples of vertex-transitive non-spherical homology spheres exist:
There are exactly three $17$-vertex triangulations $\Sigma^i_{17}$, $i=1,2,3$, 
of the Poincar\'e homology $3$-sphere $\Sigma$ that have a vertex-transitive cyclic group action.
In fact, these are the only vertex-transitive non-spherical homology $3$-spheres
with $n\leq 17$ vertices; see \cite{Lutz2003pre}.
\begin{thm}
The wreath products $\partial\Delta_d\wr \Sigma^i_{17}$
of the vertex-transitive $17$-vertex triangulations $\Sigma^i_{17}$, $i=1,2,3$,
of the Poincar\'e homology $3$-sphere $\Sigma$ give ($17d+3$)-dimensional 
vertex-transitive non-PL spheres.
\end{thm}
If instead of the wreath product we take the $k$-fold join product of these triangulations $\Sigma^i_{17}$, then
$(\Sigma^i_{17})^{*k}$ is a vertex-transitive non-PL $(4k-1)$-sphere for $k\geq 2$; cf.\ \cite{Lutz2003pre}.  In
particular, the two constructions yield examples of vertex-transitive non-PL spheres in different dimensions (unless
$d=4l$ and $k=17l+1$ for $l\geq 1$).

\subsection{Neighborly Wreath Product Spheres}

If $K$ is a simplicial $(e-1)$-sphere, then it is either the boundary of an $e$-simplex, which is $e$-neigh\-borly, or
it is at most $\lfloor\frac{e-2}{2}\rfloor$-neighborly by the van Kampen-Flores Theorem; see
Gr\"unbaum~\cite[11.1.3]{Gruenbaum1967}.  Spheres that are $\lfloor\frac{e-2}{2}\rfloor$-neighborly are simply called
\emph{neighborly}.  The wreath product $\partial\Delta_d\wr K$ of a $k$-neighborly simplicial sphere $K$ with
$\partial\Delta_d$ is a ($k(d+1)+d$)-neighborly ($nd+e-1$)-dimensional sphere by Proposition \ref{prop:neighborly}.  If
$\Delta_d$ is a point, then $\partial\Delta_d\wr K$ is neigh\-borly if and only if $K$ is neighborly, since
$\partial\Delta_d\wr K=K$.  If $K=\partial\Delta_{e+1}$, then $\partial\Delta_d\wr K$ is a simplex, which is neighborly.
\begin{prop}
  \label{prop:neighborly_up}
  Let $K$ be a simplicial $(e-1)$-sphere, different from the boundary of a simplex, and $d\geq 1$. Then
  $\partial\Delta_d\wr K$ is neighborly if and only if $K$ is neighborly and the parameters $e-1$, $n$, and $d$ obey the
  conditions that $e-1$ is odd, $e+2\leq n\leq e+3$, and $d=1$ in the case $n=e+3$.
\end{prop}
\begin{proof}
  Let $K$ be different from the boundary of a simplex (hence $n\geq e+2$), $d\geq 1$, and $\partial\Delta_d\wr K$ be
  neighborly, i.e., $\lfloor\frac{nd+e}{2}\rfloor$-neighborly.  For fixed $d$ and $e-1$, $\partial\Delta_d\wr K$ is at
  most ($\lfloor\frac{e}{2}\rfloor (d+1)+d$)-neighborly by Proposition~\ref{prop:neighborly}.  Therefore,
  $\partial\Delta_d\wr K$ can be neighborly only for small $n$.  Let $n\geq e+4$, then
\[
\begin{array}{lllll}
\displaystyle
\lfloor\frac{nd+e}{2}\rfloor 
   & \displaystyle\geq & \displaystyle\lfloor\frac{(e+4)d+e}{2}\rfloor
   & \displaystyle=    & \displaystyle\lfloor\frac{e(d+1)+4d}{2}\rfloor \\[2mm]
   & \displaystyle>    & \displaystyle\frac{e(d+1)+4d}{2}-1
   & \displaystyle\geq & \displaystyle\lfloor\frac{e}{2}\rfloor (d+1)+2d-1 \\[2mm]
   & \displaystyle\geq & \displaystyle\lfloor\frac{e}{2}\rfloor (d+1)+d.
\end{array}
\]
Thus, $\partial\Delta_d\wr K$ is not neighborly for $n\geq e+4$,
and this is also the case for $n=e+3$ and $d>1$, since then
\[
\lfloor\frac{nd+e}{2}\rfloor 
\,\, =    \,\, \lfloor\frac{(e+3)d+e}{2}\rfloor
\,\, >    \,\, \lfloor\frac{e}{2}\rfloor (d+1)+\frac{3d-2}{2}
\,\, \geq \,\, \lfloor\frac{e}{2}\rfloor (d+1)+d.
\]
For $n=e+3$ and $d=1$,
$\lfloor\frac{(e+3)\cdot 1+e}{2}\rfloor =\lfloor\frac{e}{2}\rfloor (1+1)+1$
if and only if $e-1$ is odd.
Finally, let $n=e+2$. Then 
\[
\lfloor\frac{(e+2)d+e}{2}\rfloor
\,\, = \,\,  \lfloor\frac{e(d+1)+2d}{2}\rfloor
\,\, = \,\, \lfloor\frac{e(d+1)}{2}\rfloor +d,
\]
where the last expression is equal to $\lfloor\frac{e}{2}\rfloor (d+1)+d$ 
if and only if $e-1$ is odd.
Note that if $K$ is less than $\lfloor\frac{e}{2}\rfloor$-neighborly,
then $\partial\Delta_d\wr K$ never is neighborly.
\end{proof}

\begin{cor}
All neighborly simplicial spheres that are wreath products 
$\partial\Delta_d\wr K$ of the boundary $\partial\Delta_d$ of a $d$-simplex $\Delta_d$ of dimension $d\geq 1$ 
and some simplicial sphere $K$ with $n$ vertices are polytopal.
In particular, $e+1\leq n\leq e+3$ and $K=\partial C_{e}(n)$.
If $n=e+1$, then $\partial\Delta_d\wr \partial C_{e}(n)=\partial\Delta_{nd+e}(n(d+1))$.
For odd $e-1$ and $n=e+2$ we have $\partial\Delta_d\wr \partial C_{e}(n)=\partial C_{nd+e}(n(d+1))$,
while for $n=e+3$ the wreath product $\partial\Delta_1\wr \partial C_{e}(e+3)$ is a neighborly polytopal simplicial sphere
different from $\partial C_{2e+3}(2e+6)$.
\end{cor}
\begin{proof}
Let a neighborly simplicial sphere be the wreath product
$\partial\Delta_d\wr K$ of a simplicial sphere $K$ with $n$ vertices
with the boundary of a simplex of dimension $d\geq 1$. By Proposition~\ref{prop:neighborly_up}
and the comments before, $K$ is neighborly and $e+1\leq n\leq e+3$.
According to Barnette and Gannon \cite{BarnetteGannon1976} 
every $(e-1)$-dimensional simplicial manifold with $n\leq e+4$ vertices for
$e-1=3$ and $e-1\geq 5$, and with $n\leq e+3$ vertices for $e-1=4$,
is a combinatorial sphere. Moreover, Mani \cite{Mani1972} showed that
combinatorial $(e-1)$-spheres with $n\leq e+3$ are polytopal. 
Hence, $K$ is polytopal. Is also follows from Proposition~\ref{prop:neighborly_up}
that if $K$ is not the boundary of a simplex $\Delta_{e}$ with $e+1$ vertices,
then $e-1$ is odd. Furthermore (see \cite[Ch.\ 6 \& 7]{Gruenbaum1967} for a discussion
and additional references), the number of odd-dimensional (even-dimensional) 
neighborly simplicial spheres with $n$ vertices is equal to one 
if and only $e+1\leq n\leq e+3$ ($e+1\leq n\leq e+2$).
Therefore, $K=\partial \Delta_{e}(n)$ and $\partial\Delta_d\wr \partial \Delta_{e}(n)=\partial\Delta_{nd+e}(n(d+1))$
for $n=e+1$. For odd $e-1$, $K=\partial C_{e}(n)$ for $e+2\leq n\leq e+3$
and $\partial\Delta_d\wr \partial C_{e}(n)=\partial C_{nd+e}(n(d+1))$ for $n=e+2$.
If $n=e+3$, then $d=1$ by Proposition~\ref{prop:neighborly_up}, 
so $\partial\Delta_1\wr \partial C_{e}(e+3)$ is a sphere of even dimension $2e+2$.
Since the odd-dimensional sphere $\partial C_{e}(e+3)$ has a vertex-transitive
dihedral (combinatorial and geometric) symmetry group $D_{e+3}$,
the $(2e+2)$-sphere $\partial\Delta_1\wr \partial C_{e}(e+3)$ 
with $2e+6$ vertices has the group $\ZZ_2\wr D_{e+3}$
as vertex-transitive symmetry group. However, the automorphism group
of $\partial C_{2e+3}(2e+6)$ is $\ZZ_2\times \ZZ_2$; 
cf.\ \cite{KaibelWassmer}. Thus, $\partial\Delta_1\wr \partial C_{e}(e+3)$
is distinct from $\partial C_{2e+3}(2e+6)$. 
\end{proof}

\begin{rem} The existence of this  series of odd-dimensional neighborly simplicial $(2e+3)$-polytopes $\Delta_1\wr C_{e}(e+3)$ on
  $2e+6$ vertices with a vertex-transitive symmetry group\, $\ZZ_2\wr D_{e+3}$\, for even $e\geq 2$, can also be derived
  from the results in Gr\"unbaum~\cite[\S6.2]{Gruenbaum1967}.  The numbers of different odd-dimensional neighborly
  simplicial $(2e+3)$-polytopes with $(2e+3)+3$ vertices can be found in \cite{AltshulerMcMullen1973}.  However, it
  seems to be unknown whether there are vertex-transitive neighborly simplicial polytopes other than the simplex,
  even-dimensional cyclic polytopes, and the odd-dimensional series $\Delta_1\wr C_{e}(e+3)$ for even $e\geq 2$.
  Further examples of odd-dimensional vertex-transitive neighborly simplicial spheres can be found in \cite{Lutz1999};
  for these examples it is open whether or not they are polytopal.
\end{rem}

\begin{exmp}
  For all $m\geq 1$, the cyclic polytope $C_{4m-2}(4m)$ has the following descriptions,
  $C_{4m-2}(4m)=(\Delta_{2m-1})^{*2} = \Delta_{m-1}\wr C_2(4)=\Delta_1\wr C_{2m-2}(2m)$.  In particular,
  $C_{6}(8)=(\Delta_3)^{*2}=\Delta_1\wr C_2(4)$\, and\, $C_{10}(12)=(\Delta_5)^{*2}=\Delta_2\wr C_2(4)=\Delta_1\wr
  C_4(6)$.
\end{exmp}

Note that the cyclic polytopes in odd dimensions are dual wedges over cyclic polytopes in one dimension less.

\begin{exmp}
  The $7$-polytope $\Delta_1\wr C_{2}(5)$, with $f$-vector $f=(10,45,120,205,222,140,40)$, is the smallest neighborly
  wreath product polytope which is not a cyclic polytope.  The affine Gale diagram of $\Delta_1\wr C_{2}(5)$ is
  $1$-dimensional, and it arises from the Gale diagram of the pentagon by doubling the vertices; it is displayed in
  Figure~\ref{fig:gale}.

  \begin{figure}[htbp]
    \centering
    \includegraphics[width=.5\textwidth]{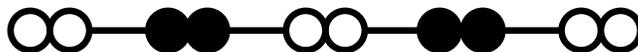}
    \caption{Affine Gale diagram of $\Delta_1\wr C_{2}(5)$.\label{fig:gale}}
  \end{figure}
\end{exmp}

\section{Recognition of Revisiting Paths}

Another interesting application of one-point suspensions is for the construction of counterexamples to the Hirsch
conjecture for simplicial spheres, which states that the diameter of the dual graph of a $(d-1)$-dimensional simplicial
sphere with $n$ vertices is bounded above by $n-d$. In fact, the original Hirsch conjecture, formulated by Hirsch in
1957 (cf.\ \cite[p.~168]{Dantzig1963}), plays an important role in the study of the computational complexity of the
simplex algorithm of linear programming (see the surveys in \cite{KleeKleinschmidt1987} and \cite{ProvanBillera1980}); it
asserts that the diameter of the graph of a $d$-polytope with $n$ facets, in other words, the number of pivot steps that
an edge-following LP algorithm needs for this polytope in the worst case with respect to a best possible choice of the
pivots, is smaller or equal to $n-d$.

While the best bound known for the diameter is super-polynomial \cite{Kalai1992,KalaiKleitman1992}, it also turns out to
be a non-trivial problem to actually construct simple (or, dually, simplicial) polytopes for all possible parameters
$(n,d)$ which can attain the Hirsch bound.  Interestingly, in the known constructions by Holt, Klee, and Fritzsche
\cite{FritzscheHolt1999,Holt2003,HoltKlee1998} (dual) wedges play a key role.

Provan and Billera \cite{ProvanBillera1980} showed that all vertex-decomposable simplicial spheres (or even more
general, all vertex-decomposable simplicial complexes) satisfy the Hirsch conjecture. Moreover, they proved that
triangulated $2$-dimensional spheres are vertex-decomposable, thus, in particular, verifying the Hirsch conjecture for
$3$-dimensional polytopes. Nevertheless, the Hirsch conjecture for (simplicial) $d$-polytopes is still open for $d\geq
4$.  For simplicial spheres the Hirsch conjecture was disproved in 1978 by Walkup \cite{Walkup1978} who provided a
$27$-dimensional counterexample with $56$ vertices.  A much smaller counterexample of dimension $11$ with $24$ vertices
was constructed by Mani and Walkup \cite{ManiWalkup1980}. Their construction is based on a $3$-dimensional sphere $D$
with $16$ vertices for which there is a pair of disjoint tetrahedra such that every path of adjacent facets joining
these two tetrahedra revisits at least one vertex that has previously been left behind.  In other words, the $3$-sphere
$D$ of Mani and Walkup provides a counterexample to the simplicial version of the $W_v$-path conjecture by Klee and
Wolfe \cite{Klee1966} ruling out such revisiting paths. If we successively one-point suspend $D$ with respect to all
vertices except for those eight vertices of the two tetrahedra for which we have the revisiting paths, then it follows
from work of Adler and Dantzig \cite{AdlerDantzig1974} (cf.\ \cite{KleeWalkup1967} and \cite{ManiWalkup1980}) that the
resulting $11$-dimensional sphere with $24$ vertices is a counterexample to the Hirsch conjecture. (In fact, if $\{
v_1,v_2,v_3,v_4\}$ and $\{ v_5,v_6,v_7,v_8\}$ are the two disjoint tetrahedra in $D$ for which we have the revisiting
paths, then $\{ v_1,v_2,v_3,v_4,v'_9\}$ and $\{ v_5,v_6,v_7,v_8,v''_9\}$ are two disjoint $4$-simplices in the one-point
suspension $\susp_1(v_9,D)$, which as well are joined by revisiting paths only, etc.)

\bibliographystyle{plain}
\bibliography{main}

\begin{thebibliography}{10}

\bibitem{AdlerDantzig1974}
I.~Adler and G.~B. Dantzig.
\newblock Maximum diameter of abstract polytopes.
\newblock {\em Math.\ Program.\ Study}, 1:20--40, 1974.

\bibitem{AltshulerMcMullen1973}
A.~Altshuler and P.~McMullen.
\newblock The number of simplicial neighbourly $d$-polytopes with $d+3$
  vertices.
\newblock {\em Mathematika\/}, 20:263--266, 1973.

\bibitem{BagchiDatta1998}
B.~Bagchi and B.~Datta.
\newblock A structure theorem for pseudomanifolds.
\newblock {\em Discrete Math.}, 188:41--60, 1998.

\bibitem{BarnetteGannon1976}
D.~Barnette and D.~Gannon.
\newblock Manifolds with few vertices.
\newblock {\em Discrete Math.}, 16:291--298, 1976.

\bibitem{Bjoerner1995}
A.~Bj\"orner.
\newblock Topological methods.
\newblock In R.~Graham, M.~Gr\"otschel, and L.~Lov\'asz, editors, {\em Handbook
  of Combinatorics}, chapter~34, pages 1819--1872. Elsevier, Amsterdam, 1995.

\bibitem{BjoernerLutz2000}
A.~Bj\"orner and F.~H. Lutz.
\newblock Simplicial manifolds, bistellar flips and a $16$-ver\-tex
  triangulation of the {Poincar\'{e}} homology $3$-sphere.
\newblock {\em Exp.\ Math.}, 9:275--289, 2000.

\bibitem{BjoernerLutz2003}
A.~Bj\"orner and F.~H. Lutz.
\newblock A $16$-vertex triangulation of the {P}oincar\'e homology $3$-sphere
  and non-{PL} spheres with few vertices.
\newblock {\em Electronic Geometry Models}, {\rm No.\ 2003.04.001}, 2003.
\newblock \url{http://www.eg-models.de/2003.04.001}.

\bibitem{Cannon1979}
J.~W. Cannon.
\newblock Shrinking cell-like decompositions of manifolds. {C}odimension three.
\newblock {\em Ann.\ Math.}, 110:83--112, 1979.

\bibitem{Chari2000}
M.~K. Chari.
\newblock On discrete {M}orse functions and combinatorial decompositions.
\newblock {\em Discrete Math.}, 217:101--113, 2000.

\bibitem{Dantzig1963}
G.~B. Dantzig.
\newblock {\em Linear Programming and Extensions}.
\newblock Princeton University Press, Princeton, NJ, 1963.

\bibitem{Edwards1975}
R.~D. Edwards.
\newblock The double suspension of a certain homology $3$-sphere is {$S^5$}.
\newblock {\em Notices AMS\/}, 22:A--334, 1975.

\bibitem{Forman1998}
R.~Forman.
\newblock Morse theory for cell complexes.
\newblock {\em Adv.\ Math.}, 134:90--145, 1998.

\bibitem{Forman2002}
R.~Forman.
\newblock A user's guide to discrete {M}orse theory.
\newblock {\em S\'emin.\ Lothar.\ Comb.}, 48:B48c, 35 p., electronic only,
  2002.

\bibitem{FritzscheHolt1999}
K.~Fritzsche and F.~B. Holt.
\newblock More polytopes meeting the conjectured {H}irsch bound.
\newblock {\em Discrete Math.}, 205:77--84, 1999.

\bibitem{Gruenbaum1967}
B.~Gr\"unbaum.
\newblock {\em Convex Polytopes}, volume~16 of {\em Pure and Applied
  Mathematics}.
\newblock Interscience Publishers, London, 1967.
\newblock Second edition (V.~Kaibel, V.~Klee, and G.~M.~Ziegler, eds.),
  Graduate Texts in Mathematics \textbf{221}. Springer-Verlag, New York, NY,
  2003.

\bibitem{HachimoriZiegler2000}
M.~Hachimori and G.~M. Ziegler.
\newblock Decompositions of simplicial balls and spheres with knots consisting
  of few edges.
\newblock {\em Math.\ Z.}, 235:159--171, 2000.

\bibitem{Holt2003}
F.~B. Holt.
\newblock Maximal nonrevisiting paths in simple polytopes.
\newblock {\em Discrete Math.}, 263:105--128, 2003.

\bibitem{HoltKlee1998}
F.~B. Holt and V.~Klee.
\newblock Many polytopes meeting the conjectured {H}irsch bound.
\newblock {\em Discrete Comput.\ Geom.}, 20:1--17, 1998.

\bibitem{KahnSaksSturtevant1984}
J.~Kahn, M.~Saks, and D.~Sturtevant.
\newblock A topological approach to evasiveness.
\newblock {\em Combinatorica\/}, 4:297--306, 1984.

\bibitem{KaibelWassmer}
V.~Kaibel and A.~Wa{\ss}mer.
\newblock Automorphism groups of cyclic polytopes.
\newblock To appear in \emph{Triangulated Manifolds with Few Vertices} by
  F.~H.~Lutz.

\bibitem{Kalai1992}
G.~Kalai.
\newblock Upper bounds for the diameter and height of graphs of convex
  polyhedra.
\newblock {\em Discrete Comput.\ Geom.}, 8:363--372, 1992.

\bibitem{KalaiKleitman1992}
G.~Kalai and D.~J. Kleitman.
\newblock A quasi-polynomial bound for the diameter of graphs of polyhedra.
\newblock {\em Bull.\ Am.\ Math.\ Soc., New Ser.}, 26:315--316, 1992.

\bibitem{Klee1966}
V.~Klee.
\newblock Paths on polyhedra. {II}.
\newblock {\em Pac.\ J.\ Math.}, 17:249--262, 1966.

\bibitem{KleeKleinschmidt1987}
V.~Klee and P.~Kleinschmidt.
\newblock The $d$-step conjecture and its relatives.
\newblock {\em Math.\ Oper.\ Res.}, 12:718--755, 1987.

\bibitem{KleeWalkup1967}
V.~Klee and D.~W. Walkup.
\newblock The $d$-step conjecture for polyhedra of dimension $d<6$.
\newblock {\em Acta Math.}, 117:53--78, 1967.

\bibitem{Lutz2003pre}
F.~H. Lutz.
\newblock \emph{Triangulated {M}anifolds with {F}ew {V}ertices}.
\newblock In preparation.

\bibitem{Lutz2003gpre}
F.~H. Lutz.
\newblock Small examples of non-constructible simplicial balls and spheres.
\newblock SIAM J. Discrete Math., to appear; \url{arXiv:math.CO/0309149}, 2003,
  9 pages.

\bibitem{Lutz1999}
F.~H. Lutz.
\newblock {\em Triangulated Manifolds with Few Vertices and
  Vertex-Tran\-si\-tive Group Actions. {\rm Dissertation}}.
\newblock Shaker Verlag, Aachen, 1999, 146 pages.

\bibitem{Lutz2002}
F.~H. Lutz.
\newblock Examples of {$\mathbb Z$}-acyclic and contractible vertex-homogeneous
  simplicial complexes.
\newblock {\em Discrete Comput.\ Geom.}, 27:137--154, 2002.

\bibitem{Lutz2004a}
F.~H. Lutz.
\newblock A vertex-minimal non-shellable simplicial $3$-ball with $9$ vertices
  and $18$ facets.
\newblock {\em Electronic Geometry Models}, {\rm No.\ 2003.05.004}, 2004.
\newblock \url{http://www.eg-models.de/2003.05.004}.

\bibitem{Lutz2004b}
F.~H. Lutz.
\newblock Vertex-minimal not vertex-decomposable balls.
\newblock {\em Electronic Geometry Models}, {\rm No.\ 2003.06.001}, 2004.
\newblock \url{http://www.eg-models.de/2003.06.001}.

\bibitem{Mani1972}
P.~Mani.
\newblock Spheres with few vertices.
\newblock {\em J.\ Comb.\ Theory, Ser.\ A\/}, 13:346--352, 1972.

\bibitem{ManiWalkup1980}
P.~Mani and D.~W. Walkup.
\newblock A $3$-sphere counterexample to the {$W_v$}-path conjecture.
\newblock {\em Math.\ Oper.\ Res.}, 5:595--598, 1980.

\bibitem{McMullen1976}
P.~McMullen.
\newblock Constructions for projectively unique polytopes.
\newblock {\em Discrete Math.}, 14:347--358, 1976.

\bibitem{Munkres1984b}
J.~R. Munkres.
\newblock Topological results in combinatorics.
\newblock {\em Mich.\ Math.\ J.}, 31:113--128, 1984.

\bibitem{Newman1926c}
M.~H.~A. Newman.
\newblock On the foundations of combinatory analysis situs. {I}, {II}.
\newblock {\em Proc.\ Royal Acad.\ Amsterdam}, 29:611--626, 627--641, 1926.

\bibitem{ProvanBillera1980}
J.~S. Provan and L.~J. Billera.
\newblock Decompositions of simplicial complexes related to diameters of convex
  polyhedra.
\newblock {\em Math.\ Oper.\ Res.}, 5:576--594, 1980.

\bibitem{ThrelfallSeifert1931}
W.~Threlfall and H.~Seifert.
\newblock Topologische {U}ntersuchung der
  {D}is\-kon\-ti\-nu\-i\-t\"ats\-be\-reiche endlicher {B}ewegungsgruppen des
  dreidimensionalen sph\"a\-risch\-en {R}aumes.
\newblock {\em Math.\ Ann.}, 104:1--70, 1931.

\bibitem{Walkup1978}
D.~W. Walkup.
\newblock The {H}irsch conjecture fails for triangulated $27$-spheres.
\newblock {\em Math.\ Oper.\ Res.}, 3:224--230, 1978.

\bibitem{WeberSeifert1933}
C.~Weber and H.~Seifert.
\newblock Die beiden {D}odekaederr\"aume.
\newblock {\em Math.\ Z.}, 37:237--253, 1933.

\bibitem{Welker1999}
V.~Welker.
\newblock Constructions preserving evasiveness and collapsibility.
\newblock {\em Discrete Math.}, 207:243--255, 1999.

\bibitem{Ziegler1995}
G.~M. Ziegler.
\newblock {\em Lectures on Polytopes}, volume 152 of {\em Graduate Texts in
  Mathematics}.
\newblock Springer-Verlag, New York, NY, 1995.
\newblock Revised edition, 1998.

\end{thebibliography}

\end{document}